\documentclass[journal,twoside]{IEEEtran} 

\usepackage{xcolor,subfigure}
\newif\ifdraft

\drafttrue
\usepackage{amssymb,amsmath, amsbsy, amstext,chemarrow,mathrsfs}
\usepackage{graphicx}
\usepackage{amsmath,stfloats}
\usepackage{tcolorbox}
\usepackage[ruled,vlined,linesnumbered]{algorithm2e}
\usepackage{bm}

\usepackage{enumitem,multirow}
\usepackage{accents}
\usepackage{etoolbox}
\makeatletter
\patchcmd{\@makecaption}
  {\scshape}
  {}
  {}
  {}
\makeatother
\usepackage{makecell}

\makeatletter
\renewcommand\normalsize{%
\@setfontsize\normalsize\@xpt\@xiipt
\abovedisplayskip 8\p@ \@plus2\p@ \@minus5\p@
\abovedisplayshortskip \z@ \@plus3\p@
\belowdisplayshortskip 8\p@ \@plus3\p@ \@minus3\p@
\belowdisplayskip \abovedisplayskip
\let\@listi\@listI}
\makeatother

\graphicspath{{./figures/}}

\newtheorem{definition}{\bfseries Definition}
\newtheorem{proposition}{\bfseries Proposition}
\newtheorem{theorem}{\bfseries Theorem}
\newtheorem{lemma}{\bfseries Lemma}

\newtheorem{remark}{\bfseries Remark}

\def\vt{{\mbox{vertex}}}

\def\diag{\operatorname{diag}}

\def\w{\bm{w}}
\def\v{\bm{v}}

\def\x{\bm{x}}

\def\u{\bm{u}}

\newcommand{\R}{\mathbb{R}}

\newcommand{\mat}[1]{\begin{bmatrix} #1 \end{bmatrix}}

\newcommand{\li}{[\![}
\newcommand{\ri}{]\!]}
\newcommand{\bbox}[1]{[\![ #1 ]\!]}

\allowdisplaybreaks[4]

\newif\ifdraft

\drafttrue

\title{\LARGE \bf
Robust Stability of Neural Network Control Systems with Interval Matrix Uncertainties}

\author{Yuhao Zhang and Xiangru Xu
\thanks{Yuhao Zhang and Xiangru Xu are with the Department of Mechanical Engineering, University of Wisconsin-Madison,
        Madison, WI, USA. Email: 
        {\tt\small \{yuhao.zhang2,xiangru.xu\}@wisc.edu}.}%
}

\begin{document}
\maketitle
\begin{abstract}
Neural networks have become increasingly popular in controller design due to their versatility and efficiency. However, their integration into feedback systems can pose stability challenges, particularly in the presence of uncertainties. This work addresses the problem of certifying robust stability in neural network control systems with interval matrix uncertainties. Leveraging classical robust stability techniques and the quadratic constraint-based method to characterize the input-output behavior of neural networks, we derive novel robust stability certificates formulated as linear matrix inequalities. To reduce computational complexity, we introduce three relaxed sufficient conditions and establish their equivalence in terms of feasibility. Additionally, we explore their connections to existing robust stability results. The effectiveness of the proposed approach is demonstrated through inverted pendulum and mass-spring-damper examples.
\end{abstract}

\section{Introduction}\label{sec:intro}

Neural Networks (NNs) are widely used in controller design for various dynamical systems because of their universal approximation capabilities {\cite{fazlyab2020safety,hornik1989multilayer,nikolakopoulou2022dynamic,schwan2023stability,zhang2022safety,zhang2023backward,zhang2024reachability}. However, Neural Network Control Systems (NNCSs), which are feedback systems with NN controllers in the loop, often lack formal stability guarantees due to the complex nature of NNs. This issue is even more critical when dealing with uncertainties in system models, as NNs are known to be sensitive to perturbations \cite{goodfellow2014explaining}. As a result, it is essential to certify the properties of NNCSs such as stability and safety before deploying them in practical applications.

Recently, there has been a growing interest in addressing the stability verification problem of NNCSs. Various methods have emerged to certify the stability by constructing candidate Lyapunov functions, either through optimization techniques \cite{fabiani2022reliably,karg2020stability,newton2023sparse,schwan2023stability} or self-supervised learning approaches \cite{dawson2023safe}. In recent works such as \cite{de2023event,fazlyab2020safety,hu2020reach,jin2020stability,pauli2021linear,yin2021stability,yin2021imitation}, Quadratic Constraints (QCs) were utilized in analyzing NNs by abstracting nonlinear activation functions. This approach enables the formulation of safety verification and robustness analysis for NNs against norm-bounded perturbations as semi-definite programs \cite{fazlyab2020safety}. For NNCSs without uncertainties, QCs were also used for forward reachability analysis \cite{hu2020reach} and stability analysis \cite{jin2020stability,yin2021imitation}. When uncertainties exist in the system dynamics, the QC-based methodology was further applied to the robust stability analysis of NNCSs incorporating perturbations that align with Integral Quadratic Constraints (IQCs) \cite{de2023event,yin2021stability}. Along this line, an improved stability analysis via acausal Zames-Falb multipliers was considered in \cite{pauli2021linear}, which offers enhanced stability assurances and the potential for larger Regions Of Attraction (ROAs). These existing results are important but they all require the uncertainties to be represented by IQCs, which limits their application to other commonly used uncertain system models.

\begin{figure}[t]
    \centering
    \includegraphics[width=\linewidth]{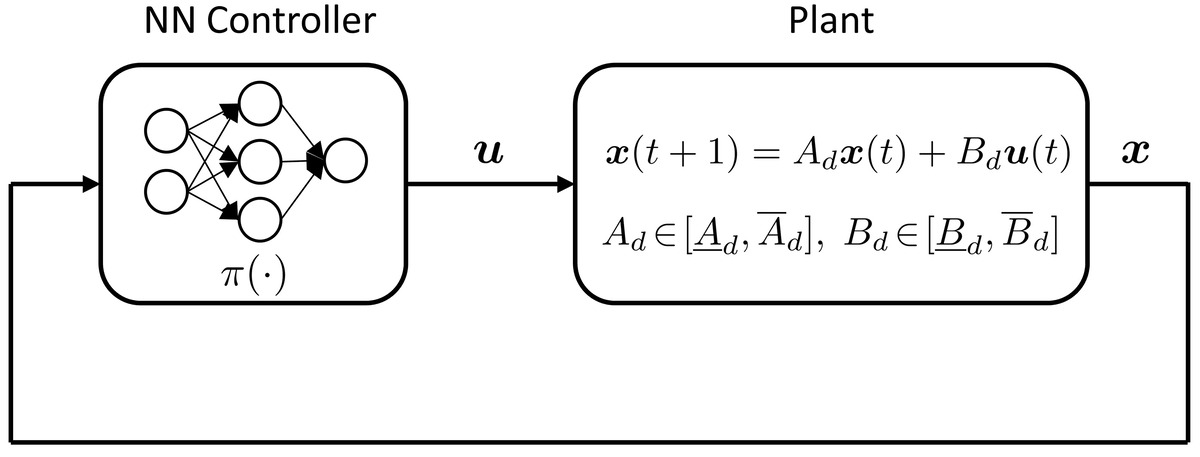}
    \caption{Feedback system that contains a plant with interval matrix uncertainties and NN controller $\pi$.}
    \label{fig:flow}
\end{figure}

In this work, we consider the problem of verifying the robust stability of NNCSs with interval matrix uncertainties (see  Fig. \ref{fig:flow}). 
For conventional systems without NN components, there exists a substantial body of work addressing robust stability in the presence of uncertainties, such as \cite{petersen1987stabilization,petersen1986riccati}. The readers may refer to \cite{petersen2014robust} and the references therein, for additional pointers to the relevant literature. In the case of linear systems with interval matrix uncertainties, various necessary and sufficient conditions were developed in the form of Linear Matrix Inequalities (LMIs) \cite{alamo2008new,ben2002tractable,mao2003quadratic}. Although \cite{alamo2008new} offers a concise comparison of the performance of different sufficient conditions, a clear connection among the LMI conditions remains to be fully explored. 
Leveraging insights from these classic robust stability results for systems without NN components and the recent advances in QC-based techniques for NN controllers, we introduce novel and efficient methods to certify the robust stability of uncertain NNCSs.

The contributions of this work are at least threefold: (i) By using the QC-based abstraction techniques and the vertices of the interval matrix uncertainties, a certificate is formulated as a finite number of LMIs to verify stability and approximate the ROAs of uncertain NNCSs; (ii) To reduce the computation burden associated with solving the LMIs, three relaxed sufficient LMI conditions are proposed featuring fewer numbers of decision variables and smaller sizes of LMI constraints; (iii) The equivalence of feasibility among the three relaxed LMIs is established, and their connections with existing robust stability results for interval matrix uncertainties are also built. 
The remainder of this paper is laid out as follows: Section \ref{Prelim} presents the problem formulation and preliminaries on stability analysis for nominal NNCSs. Section \ref{sec:uncertain} introduces a sufficient stability condition for uncertain NNCSs in the form of LMIs, which is then relaxed into three LMIs with fewer computation complexities in Section \ref{sec:relax}. Section \ref{sec:exp} features two numerical simulations, followed by concluding remarks in Section \ref{sec:concl}.

\textit{Notation:} 
The $i$-th entry of a vector $\x\in \mathbb{R}^n$ is denoted by $x_i$ with $i\in \bbox{n}$, where $ \bbox{n}\triangleq \{1,\dots,n\}$. 
For a matrix $A\in \mathbb{R}^{n\times m}$, $A(i,j)$ denotes the $i$-th row and $j$-th column entry of $A$. 
Given a square matrix $A\in \mathbb{R}^{n\times n}$, we denote the $(i, j)$ cofactor of $A$ as $C^A_{i,j}$ and the minimum eigenvalue of $A$ as $\lambda_{min}(A)$. The $n\times n$ identity matrix is denoted as ${I}_{n}$ and $\bm{e}_i^n$ is the $i$-th column of ${I}_{n}$. The $n\times m$ matrix whose entries are all $0$ (resp., $1$) is denoted as $\bm{0}_{n\times m}$ (resp., $\bm{1}_{n\times m}$); we may omit the subscripts when the dimensions are evident from the context.
The sets of symmetric matrices, positive semi-definite matrices, and positive definite matrices are denoted as $\mathbb{S}^{ n}$, $\mathbb{S}^{ n}_{\succeq 0}$, and $\mathbb{S}^{n}_{\succ 0}$, respectively. The set of non-negative vectors is denoted as $\mathbb{R}^{n}_{\geq 0}$. 
The notation $\geq$ and $\leq$ are used to denote the entry-wise relationship of matrices and vectors with appropriate dimensions.
Given matrices $\underline{A},\overline{A} \in \mathbb{R}^{n\times m}$ with $\underline{A}(i,j)\leq\overline{A}(i,j)$ for any $i\in \bbox{n}$ and $j \in \bbox{m}$, the interval matrix $[\underline{A},\overline{A}]$ is defined as $[\underline{A},\overline{A}] \triangleq \{A \in \mathbb{R}^{n\times m}| \underline{A}(i,j)\leq A(i,j) \leq \overline{A}(i,j),\forall i\in \bbox{n}, \forall j\in \bbox{m}\}$. The set of vertices of an interval matrix is defined as $\vt([\underline{A},\overline{A}]) \triangleq \{A \in \mathbb{R}^{n\times m}| A(i,j)=\underline{A}(i,j) \text{ or } \overline{A}(i,j), \forall i\in \bbox{n}, \forall j\in \bbox{m}\}$. 
The symbol $*$ denotes entries whose values follow from symmetry.

\section{Preliminary}\label{Prelim}

\subsection{Problem statement}

In this paper, we consider the following discrete-time linear system with interval matrix uncertainties:
\begin{equation}\label{dt-sys}
    \x{(t+1)} =  A_d \x(t) +  B_d \u(t)
\end{equation}
where $\x(t)\in  \mathbb{R}^n,\; \u(t)\in \mathbb{R}^m$ denote the state and the control input, respectively. The system matrices $A_d,B_d$ are subject to interval matrix uncertainties with known lower and upper bounds:
\begin{align}\label{equ:uncertain}
A_d \in [\underline{A}_d,\overline{A}_d] \subset\mathbb{R}^{n\times n},\;B_d \in[\underline{B}_d,\overline{B}_d] \subset\mathbb{R}^{n\times m}.
\end{align}
The system described by \eqref{dt-sys} with interval matrix uncertainties \eqref{equ:uncertain} can be utilized in many practical applications where uncertainties arise from unpredictable factors such as measurement errors, environmental disturbances, and variability in plant characteristics. Examples include robust aircraft control under varying operating conditions with uncertain aerodynamic coefficients, automotive control systems coping with tire parameter variations due to wear and tear, and robotic manipulators dealing with uncertainties in joint stiffness and damping coefficients. 

Let $A_0 = \frac{1}{2}(\overline{A}_d+\underline{A}_d)$, $A_r =  \frac{1}{2}(\overline{A}_d-\underline{A}_d)$ and $B_0 = \frac{1}{2}(\overline{B}_d+\underline{B}_d)$, $B_r =  \frac{1}{2}(\overline{B}_d-\underline{B}_d)$. Then, similar to \cite{mao2003quadratic}, the matrices $A_d$ and $B_d$ can be expressed as 
\begin{align}
& A_d = A_0 +\sum_{i,j=1}^n \bm{e}_i^n f_{ij} {(\bm{e}_j^n)}^{\top},\\
& B_d = B_0 + \sum_{i=1}^n\sum_{j=1}^m \bm{e}_i^n g_{ij}  {(\bm{e}_j^m)}^{\top},
\end{align}
with $ f_{ij}\in\mathbb{R}, |f_{ij}|\leq A_r (i,j)$ and $g_{ij}\in\mathbb{R}, |g_{ij}|\leq B_r (i,j)$. 
The controller $\u(t)$ is given as 
\begin{equation}
\u(t) = \pi(\x(t))\label{dt-input}    
\end{equation}
where $\pi: \mathbb{R}^{n} \rightarrow \mathbb{R}^{m}$  is a known $\ell$-layer Feedforward Neural Network (FNN) defined as follows: 
\begin{equation}\label{NN}
    \begin{aligned}
\bm{w}^{0}(t) &=\x(t), \\
\bm{w}^{i}(t) &=\phi^i\left(W^{i} \bm{w}^{i-1}(t)+\bm{b}^{i}\right),\ \forall i\in \bbox{\ell}, \\
\u(t) &=W^{\ell+1} \bm{w}^{\ell}(t)+\bm{b}^{\ell+1}.
\end{aligned}
\end{equation}
Here $\bm{w}^i \in \mathbb{R}^{n_i}\; (i\in\bbox{\ell})$ is the output (activation) from the $i^{th}$ layer. 
For each layer, the operations are defined by a weight matrix $W^i \in \mathbb{R}^{n_i \times n_{i-1}}$, a bias vector $\bm{b}^i \in \mathbb{R}^{n_i}$, and an activation function $\phi^i: \mathbb{R}^{n_i} \rightarrow \mathbb{R}^{n_i}$  given as
\begin{equation}
    \phi^{i}(\bm{v})\triangleq\left[\varphi\left(v_{1}\right), \cdots, \varphi\left(v_{n_{i}}\right)\right]^{\top} \label{activationphi}
\end{equation}
where $\varphi: \mathbb{R} \rightarrow \mathbb{R}$ is a scalar activation function (e.g., ReLU, sigmoid, tanh, leaky ReLU). 
Although in this work we assume  $\varphi$ is identical in all layers, the results can be easily extended to the case where different activation functions exist across different layers with minor notation changes, as long as all involved activation functions satisfy the sector-bounded properties and can be abstracted using QCs as detailed in the next subsection.

The uncertain NNCS consisting of system \eqref{dt-sys} and controller \eqref{dt-input} is a closed-loop system denoted as
\begin{equation}\label{close-sys}
\begin{aligned}
    \x(t+&1) =  A_d\x(t) +  B_d \pi(\x(t)) \\
   \mbox{where}\quad &A_d \in [\underline{A}_d,\overline{A}_d],\; B_d \in[\underline{B}_d,\overline{B}_d].
\end{aligned}
\end{equation}
We assume $\x\in\mathcal{X}$ where $\mathcal{X}\subset \mathbb{R}^n$ is called the state set. We also assume the FNN satisfies  $\pi(\bm 0) = \bm 0$, which ensures that $\x_* = \bm 0$ is an equilibrium point of the uncertain NNCS \eqref{close-sys}, i.e., 
\begin{eqnarray}\label{equ:equil}
    \x_* \!=\! A_d \x_*\!+\!B_d \pi(\x_*), \forall A_d \!\in\! [\underline{A}_d,\overline{A}_d], \forall B_d\!\in\![\underline{B}_d,\overline{B}_d].
\end{eqnarray}
Given an initial state $\x_0$, $\mathscr{X}(t,x_0,A_d,B_d)$ denotes the solution of the uncertain NNCS \eqref{close-sys} at time $t$ with $A_d\in[\underline{A}_d, \overline{A}_d]$ and $B_d\in[\underline{B}_d, \overline{B}_d]$. 

In this work, we aim to solve the following problem: \emph{Given the uncertain NNCS \eqref{close-sys} comprised of the discrete-time linear system \eqref{dt-sys} with interval matrix uncertainties \eqref{equ:uncertain} and the controller \eqref{dt-input} represented by an $\ell$-layer FNN \eqref{NN}, develop efficient conditions for verifying its local robust stability around the equilibrium point at the origin.}

\subsection{Stability of nominal NNCS}
In this section, we review the stability result for nominal NNCSs given in \cite{yin2021stability}. 
It is useful to isolate the nonlinear activation functions from the linear operations of the NN as done in \cite{fazlyab2020safety}, \cite{kim2018standard}, and \cite{yin2021stability}. Define $\bm{v}^i$ as the input to the activation function $\phi^i$:
\begin{equation}
\bm{v}^{i}(t)\triangleq W^{i} \bm{w}^{i-1}(t) + \bm{b}^i,\ \forall i\in \bbox{\ell}.
\end{equation}
Define the concatenated forms of the inputs $\v_{\phi}\in \mathbb{R}^{n_{\phi}}$ and outputs $\w_{\phi}\in \mathbb{R}^{n_{\phi}}$ of all activation functions, as well as the combined nonlinearity $\phi: \R^{n_\phi}\rightarrow \R^{n_\phi}$, respectively, as
\begin{equation}\label{vwphi}
\v_{\phi}\triangleq\mat{
\v^{1} \\
\vdots \\
\v^{\ell}
}, \;\w_{\phi}\triangleq\mat{
\w^{1} \\
\vdots \\
\w^{\ell}},\;\phi\left(\v_{\phi}\right)\triangleq\mat{
\phi^{1}\left(\v^{1}\right) \\
\vdots \\
\phi^{\ell}\left(\v^{\ell}\right)}
\end{equation}
where $n_{\phi}\triangleq\sum\limits_{i=1}^{\ell} n_{i}$. 
The scalar activation function $\varphi$ as shown in \eqref{activationphi} is applied element-wise to each entry of $\v_{\phi}$. Then the output $\w_{\phi}(t)$ can be expressed compactly as 
\begin{align}\label{wphi}
\w_{\phi}(t)=\phi\left(\v_{\phi}(t)\right).    
\end{align}
The controller $\pi$ defined in \eqref{NN} and $\v_{\phi},\w_{\phi}$ defined in \eqref{vwphi}-\eqref{wphi} can be rewritten as
\begin{equation}
\mat{\v_{\phi}(t)\\\u(t)}=N\mat{\x(t) \\ \w_{\phi}(t) \\ 1} \\ 
\end{equation}
where 
\begin{equation*}
\begin{aligned}
    N& =\left[\begin{array}{c|cccc|c}
 W^{1} & \bm 0 & \cdots & \bm 0 & \bm 0 & \bm b^{1}\\
\bm 0 & W^{2} & \cdots & \bm 0 & \bm 0 & \bm b^{2}\\
\vdots & \vdots & \ddots & \vdots & \vdots & \vdots\\
\bm 0 & \bm 0 & \cdots & W^{\ell} & \bm 0 & \bm b^{\ell}\\
\hline \bm 0 & \bm 0 & \bm 0 & \cdots & W^{\ell+1} & \bm b^{\ell+1} \\
\end{array}\right]\\ 
&\triangleq \left[\begin{array}{c|c|c}
N_{v x} & N_{v w} & N_{vb}\\
\hline N_{u x} & N_{u w} & N_{ub}\\
\end{array}\right].
\end{aligned}
\end{equation*}
The equilibrium point $\x_*= \bm 0$ can be propagated through the FNN (6) to obtain equilibrium values $\bm v_*^i, \bm w_*^i\; (i\in\bbox{\ell})$ for the inputs and outputs of each activation function, yielding $\bm v_\phi = \bm v_*$, $\bm w_\phi = \bm w_*$ 
and $\bm u_* = \bm 0$. Then, by construction, $(\x_*,\u_*,\bm v_*, \bm w_*)$ are unique and 
satisfy the following conditions for all $A_d \in [\underline{A}_d,\overline{A}_d]$,  $B_d\in[\underline{B}_d,\overline{B}_d]$:
    \begin{eqnarray}\label{equ:NNequ}
    \x_* &\!=\! A_d \x_*\!+\!  B_d \u_*,\;  
\mat{\bm v_*\\
\bm u_*
}  \!=\!N\mat{
\bm x_* \\
\bm w_* \\
1
}, \;
\bm w_*  \!=\!\phi\left(\bm v_*\right) .
\end{eqnarray}
Although system (8) may have equilibrium points other than $\x_* = \bm 0$, we focus on its robust stability around the origin $\x_* = \bm 0$ in this work.

The nonlinearities of the FNN $\pi$ imposed by the activation function can be abstracted using QCs \cite{fazlyab2020safety,megretski1997system,xu2020observer}. A list of QCs encoding various properties of activation functions can be found in \cite{fazlyab2020safety}. To simplify the stability analysis of the uncertain NNCS \eqref{close-sys}, we use QCs to abstract the sector-bounded properties of the activation function around the equilibrium point.

\begin{definition} \cite{fazlyab2020safety,yin2021stability}
Let $\alpha, \beta, \underline{v}, \overline{v}, v_* \in \mathbb{R}$ be given with $\alpha \leq \beta$ and $\underline{v} \leq v_* \leq \overline{v}$. The function $\varphi: \mathbb{R} \rightarrow \mathbb{R}$ is locally sector-bounded in the sector $[\alpha, \beta]$ around the point $\left(v_*, \varphi\left(v_*\right)\right)$ if 
$
(\Delta \varphi(v)-\alpha \Delta v) \cdot(\Delta \varphi(v)-\beta \Delta v) \leq 0, \forall v \in[\underline{v}, \overline{v}],
$ 
where $\Delta \varphi(v)\triangleq\varphi(v)-\varphi\left(v_*\right)$ and $\Delta v\triangleq v-v_*$.
\end{definition}

The sectors can be stacked into vectors $\bm{\alpha}_\phi, \bm{\beta}_\phi \in \mathbb{R}^{n_\phi}$ that provide QCs satisfied by the combined nonlinearity $\phi$.

\begin{lemma}\cite[Lemma 1]{yin2021stability} \label{lem:qc}
Let $\bm \alpha_\phi, \bm \beta_\phi, \underline{\bm v}, \overline{\bm v}, \bm v_* \in \mathbb{R}^{n_\phi}$ be given with $\bm \alpha_\phi \leq \bm \beta_\phi$, $ \underline{\bm v} \leq \bm v_* \leq  \overline{\bm v}$, and $\bm w_*\triangleq\phi\left(\bm v_*\right)$. Assume $\phi$ is locally sector-bounded in the sector $\left[\bm \alpha_\phi, \bm \beta_\phi\right]$ around the point $\left(\bm v_*, \bm w_*\right)$ entry-wise for all $\bm v_\phi \in[ \underline{\bm v},  \overline{\bm v}]$. For any  $\bm \lambda \in \mathbb{R}^{n_\phi}_{\geq 0}$ and $\bm  v_\phi \in[\underline{\bm v}, \overline{\bm v}]$, it holds that
\begin{align*}
\mat{
\bm v_\phi-\bm v_* \\
\bm w_\phi-\bm w_*
}^{\top} \Psi_\phi^{\top} M_\phi(\bm \lambda) \Psi_\phi\mat{
\bm v_\phi-\bm v_* \\
\bm w_\phi-\bm w_*
} \geq 0  
\end{align*}
where $\bm w_\phi=\phi\left(\bm v_\phi\right)$ and 
\begin{align*}
\!\!\!\!\Psi_\phi\!=\!\mat{
\!\operatorname{diag}\left(\bm \beta_\phi\right) & -I_{n_\phi} \\
\!-\!\operatorname{diag}\left(\bm \alpha_\phi\right) & I_{n_\phi}
},
M_\phi(\bm \lambda)\!=\!\mat{
\bm 0_{n_\phi\times n_\phi} & \operatorname{diag}(\bm \lambda) \\
\operatorname{diag}(\bm \lambda) & \bm 0_{n_\phi\times n_\phi}
}\!\!.
\end{align*}
\end{lemma}

In order to apply Lemma \ref{lem:qc}, we assume the bounds $\underline{\bm v},\overline{\bm v}\in\mathbb{R}^{n_\phi}$ on the activation input $\bm v_\phi$ are given (see, e.g., \cite{gowal2019scalable}). 

The following lemma provides a sufficient stability condition for the nominal NNCS \eqref{close-sys} without uncertainties.

\begin{lemma}\label{lem:lya}\cite[Theorem 1]{yin2021stability}
Consider the nominal NNCS \eqref{close-sys} without uncertainties, i.e., $A_d = \underline{A}_d = \overline{A}_d$ and $B_d = \underline{B}_d = \overline{B}_d$. Let the equilibrium point $\bm x_*=\bm 0$ and $(\bm x_*,\u_*,\bm v_*, \bm w_*)$ satisfy \eqref{equ:NNequ}. Let $\overline{\bm v}^1 \in \mathbb{R}^{n_1}, \underline{\bm v}^1\triangleq2 \bm v_*^1-\overline{\bm v}^1$, and let $\bm \alpha_\phi, \bm \beta_\phi \in \mathbb{R}^{n_\phi}$ be given vectors such that the combined nonlinearity $\phi$ is locally sector-bounded in the sector $\left[\bm \alpha_\phi, \bm \beta_\phi\right]$ around the point $\left(\bm v_*, \bm w_*\right)$. If there exists $P \in \mathbb{S}_{\succ 0}^{ n}$ and $\bm \lambda \in \mathbb{R}^{n_\phi}_{\geq 0}$ 
such that 
\begin{align}
& R_V^{\top}\left[\begin{array}{cc}
A_d^{\top} P A_d-P & A_d^{\top} P B_d \\
* & B_d^{\top} P B_d
\end{array}\right] R_V +
X(\bm \lambda)\prec 0, \label{equ:nom1}\\ 
& {\left[\begin{array}{cc}
\left(\overline{v}_i^1-v_{*, i}^1\right)^2 & W_i^1 \\
* & P
\end{array}\right] \succeq 0, \quad \forall i\in \bbox{ n_1}},\label{equ:nom2}
\end{align}
where $W_i^1$ is the $i$-th row of $W^1$, $v_{*, i}^1$ is the $i$-th entry of $\bm v_{*}^1$, and 
\begin{align}
R_V&=\mat{
I_{n} & \bm 0_{n \times n_\phi} \\
N_{u x} & N_{u w}},\; R_\phi=\mat{
N_{v x} & N_{v w} \\
\bm 0_{n_\phi \times n} & I_{n_\phi}},\label{rvrphi}\\
X(\bm \lambda) &= R_\phi^{\top} \Psi_\phi^{\top} M_\phi(\bm \lambda) \Psi_\phi R_\phi, \label{eqxlambda}   
\end{align}
then, (i) the nominal NNCS is locally stable around $\x_*=\bm 0$, and (ii) the set 
$\mathcal{E}\left(P, \x_*\right)\triangleq \{\x\in \mathbb{R}^{n}|(\x-\x_*)^{\top} P (\x-\x_*)\leq 1\}$  
is an inner-approximation of the ROA $\mathcal{R}$, i.e., $\mathcal{E}\left(P, \x_*\right) \subseteq \mathcal{R} \triangleq \{\x\in \mathbb{R}^n | \lim_{t\rightarrow \infty}\mathscr{X}(t,\x,A_d,B_d)=\x_*, A_d = \underline{A}_d = \overline{A}_d,B_d = \underline{B}_d = \overline{B}_d\}$.
\end{lemma}

\section{Stability of uncertain neural network control systems}\label{sec:uncertain}

Although Lemma \ref{lem:lya} provides a sufficient stability condition for the NNCS \eqref{close-sys} when the uncertainties are not present, it cannot be directly extended to the case with interval matrix uncertainties because of the quadratic terms involving $A_d$ and $B_d$ in \eqref{equ:nom1}. 

In this section, we first introduce an LMI condition in Proposition \ref{prop1} that linearly depends on the system matrices $A_d$ and $B_d$ for the nominal NNCS. Then, in Theorem \ref{thm:vertex}, we extend this result to the uncertain NNCS \eqref{close-sys} by formulating a sufficient robust stability certification using vertex LMI conditions. To enhance computational efficiency, three relaxed LMI conditions are presented in Section \ref{sec:relax}, which are shown to be sufficient robust stability certificates in Theorem \ref{thm:relax1} with the same level of conservativeness as proved in Proposition \ref{prop:lmi_eq}.

\begin{proposition}\label{prop1}
Let $A_d,B_d$ be two given matrices without uncertainties (i.e., $A_d = \underline{A}_d= \overline{A}_d$ and $B_d = \underline{B}_d= \overline{B}_d$), and let $R_V$, $X(\bm \lambda)$ be defined as in \eqref{rvrphi} and \eqref{eqxlambda}, respectively. 
For any $P \in \mathbb{S}_{\succ 0}^{ n}$ and $\bm \lambda \in \mathbb{R}^{n_\phi}_{\geq 0}$, the following LMI \eqref{equ:proplmi3} is equivalent to \eqref{equ:nom1}:
\begin{equation}\label{equ:proplmi3}
\mat{\!-R_V^{\top} \!\mat{I_n \\\bm 0_{m\times n}}\! P [I_n \;\; \bm 0_{n\times m}]R_V + X(\bm \lambda) & * \\ P[A_d \;\; B_d]R_V & -P} \!\prec\! 0.
\end{equation}
\end{proposition}
\begin{IEEEproof}
The left-hand side of \eqref{equ:nom1} can be written equivalently as 
\begin{align*}
   & R_V^{\top} \mat{A_d^{\top} PA_d-P & A_d^{\top}PB_d\\ * & B_d^{\top}PB_d}R_V+X(\bm \lambda) \\
    = & R_V^{\top} \mat{A_d^{\top} PA_d & A_d^{\top}PB_d\\ * & B_d^{\top}PB_d}R_V - R_V^{\top} \mat{P & \bm 0 \\\bm 0 & \bm 0}R_V\!+\!X(\bm \lambda)\\
    = & R_V^{\top} \mat{A_d^{\top} \\B_d^{\top}} P [A_d \; B_d]R_V \!-\! R_V^{\top} \mat{I_n \\\bm 0} P [I_n \;\; \bm 0]R_V \!+\!X(\bm \lambda).
\end{align*}
Since $-P\prec 0$, it is easy to check that LMI \eqref{equ:proplmi3} is equivalent to \eqref{equ:nom1} by using the Schur complement of the last equation above.
\end{IEEEproof}

Since the LMI condition \eqref{equ:proplmi3} in Proposition \ref{prop1} is linear in $A_d$ and $B_d$, it can be extended to the uncertain NNCS \eqref{close-sys}. Specifically, consider the uncertain NNCS \eqref{close-sys} with $A_d, B_d$ satisfying the interval matrix uncertainty constraints \eqref{equ:uncertain} and let $R_V$, $X(\bm \lambda)$ be defined as above. If there exists a matrix $P \in \mathbb{S}_{\succ 0}^{ n}$ and a vector $\bm \lambda \in \mathbb{R}^{n_\phi}_{\geq 0}$ such that \eqref{equ:proplmi3} holds for  any $A_d \in [\underline{A}_d,\overline{A}_d]$ and any $B_d \in [\underline{B}_d,\overline{B}_d]$, then the uncertain NNCS \eqref{close-sys} is locally stable around $\x_*=\bm 0$. This 
can be easily proven by showing that the Lyapunov function $\mathcal{V}(\x)\triangleq (\x-\x_*)^{\top} P (\x-\x_*)$ satisfies $\mathcal{V}(\x(t+1)) < \mathcal{V}(\x(t))$, for any $\x(t) \neq \x_*$, $A_d \in [\underline{A}_d,\overline{A}_d]$, $B_d \in [\underline{B}_d,\overline{B}_d]$, and $t \geq 0$. However, verifying uncertain LMIs \eqref{equ:proplmi3} for  any $A_d \in [\underline{A}_d,\overline{A}_d]$ and any $B_d \in [\underline{B}_d,\overline{B}_d]$ is known to be  NP-hard \cite{nemirovskii1993several}.

The following theorem proposes a method for solving the uncertain LMIs by replacing the interval matrix uncertainties with the vertex matrices.

\begin{theorem}\label{thm:vertex}
Consider the uncertain NNCS \eqref{close-sys} with $A_d, B_d$ satisfying the interval matrix uncertainty constraints \eqref{equ:uncertain}. Let the equilibrium point $\bm x_*=\bm 0$ and $(\bm x_*,\u_*,\bm v_*, \bm w_*)$ satisfy \eqref{equ:NNequ}, and let $\overline{\bm v}^1 \in \mathbb{R}^{n_1}$. Let $R_V$, $X(\bm \lambda)$ be defined as in \eqref{rvrphi} and \eqref{eqxlambda}, respectively.
If there exists a matrix $P \in \mathbb{S}_{\succ 0}^{ n}$ and a vector $\bm\lambda \in \mathbb{R}^{n_\phi}_{\geq 0}$, such that \eqref{equ:nom2} holds and \eqref{equ:proplmi3} holds for any $A_d \in \vt([\underline{A}_d,\overline{A}_d])$ and $B_d \in \vt([\underline{B}_d,\overline{B}_d])$, 
then, (i) the uncertain NNCS \eqref{close-sys} is locally stable around $\x_*=\bm 0$, and (ii) the set $\mathcal{E}\left(P, \x_*\right)\triangleq \{\x\in \mathbb{R}^{n}|(\x-\x_*)^{\top} P (\x-\x_*)\leq 1\}$ 
is an inner-approximation of the robust ROA $\tilde{\mathcal{R}}$, i.e., $\mathcal{E}\left(P, \x_*\right) \subseteq \tilde{\mathcal{R}} \triangleq \{\x\in \mathbb{R}^n | \lim_{t\rightarrow \infty}\mathscr{X}(t,\x,A_d,B_d)=\x_*, \forall A_d \in [\underline{A}_d , \overline{A}_d], \forall B_d \in [\underline{B}_d,\overline{B}_d]\}$. 
\end{theorem}

\begin{IEEEproof}
Denote $\vt([\underline{A}_d,\overline{A}_d]) \triangleq \{A_1,\dots,A_{n_A}\}$ and $\vt([\underline{B}_d,\overline{B}_d]) \triangleq \{B_1,\dots,B_{n_B}\}$. 

Since $[\underline{A}_d,\overline{A}_d]$ and $[\underline{B}_d,\overline{B}_d]$ are convex sets, for arbitrary  $A_d \in [\underline{A}_d,\overline{A}_d]$ and $B_d \in [\underline{B}_d,\overline{B}_d]$, there always exist non-negative scalars $\alpha_1,\dots,\alpha_{n_A},\beta_1,\dots,\beta_{n_B}$ such that $\sum_{i=1}^{n_A}\alpha_i = 1$, $\sum_{j=1}^{n_B}\beta_j = 1$, $A_d = \sum_{i=1}^{n_A}\alpha_i A_i$ and $B_d = \sum_{j=1}^{n_B}\beta_j B_j$.
It is easy to check that 
\begin{align*}
    P[A_d\;B_d]R_V = &   P[\sum_{i=1}^{n_A}\alpha_i A_i\; \sum_{j=1}^{n_B}\beta_j B_j]R_V \\
    = &  \sum_{i=1}^{n_A}\sum_{j=1}^{n_B} \alpha_i\beta_j P[A_i\;B_j]R_V.
\end{align*}
Since \eqref{equ:proplmi3} holds for any $A_d\in\{A_i\}_{i=1}^{n_A}$ and $B_d\in\{B_j\}_{j=1}^{n_B}$, it is easy to see that the left-hand side of \eqref{equ:proplmi3} is a summation of $n_A\cdot n_B$ negative definite matrices, which implies that \eqref{equ:proplmi3} holds for any $A_d \in [\underline{A}_d,\overline{A}_d]$ and any $B_d \in [\underline{B}_d,\overline{B}_d]$. 

Using Proposition \ref{prop1}, we know \eqref{equ:nom1} also holds for any $A_d \in [\underline{A}_d,\overline{A}_d]$ and any $B_d \in [\underline{B}_d,\overline{B}_d]$. 

Then, following the proof of  \cite[Theorem 1]{yin2021stability}, we can show that the Lyapunov function $\mathcal{V}(\x)\triangleq (\x-\x_*)^{\top} P (\x-\x_*)$ satisfies $\mathcal{V}(\x(t+1)) < \mathcal{V}(\x(t))$, for any $\x(t)\neq \x_*$, $A_d \in [\underline{A}_d,\overline{A}_d]$, $B_d \in [\underline{B}_d,\overline{B}_d]$, and $t \geq 0$. This completes the proof.
\end{IEEEproof}

In the following sections, we will denote the LMI conditions established in Theorem \ref{thm:vertex} as $\mbox{\textup{(LMI-Vertex)}}$.  Theorem \ref{thm:vertex} provides a sufficient condition based on the vertices of the interval uncertain matrices to certify the robust stability of the uncertain NNCS \eqref{close-sys}. While this vertex-based approach is manageable for uncertain NNCSs with relatively low state and input dimensions, in the worst-case scenario, it demands satisfaction of $2^{n(n+m)}$ vertex constraints. Although techniques such as vertex reduction methods presented in \cite{alamo2008new,calafiore2008reduced} can alleviate the computational load, the number of LMI constraints still exhibits exponential growth rates with respect to the system dimension $n$ and input dimension $m$.

\section{Relaxed sufficient stability conditions}\label{sec:relax}

To avoid the vertex enumeration involved in Theorem \ref{thm:vertex}, this section will introduce three relaxed sufficient conditions for certifying the robust stability of the uncertain NNCS \eqref{close-sys}. 
Since intervals are closed under multiplication \cite{jaulin2001interval}, we can compute the multiplication 
$
    [\underline{B}_d,\overline{B}_d] N_{uw} = [\operatorname{inf}_{B_d \in [\underline{B}_d,\overline{B}_d]}  B_d N_{uw},\; \operatorname{sup}_{B_d \in [\underline{B}_d,\overline{B}_d]} \allowbreak  B_d \allowbreak N_{uw}] \triangleq [\underline{{B}}_N,{\overline{B}}_N].
$
Let 
\begin{align}\label{equ:B_new_0}
    \tilde{B}_0 = \frac{1}{2}(\overline{{B}}_N+\underline{{B}}_N) \text{, and } \tilde{B}_r = \frac{1}{2}(\overline{{B}}_N-\underline{{B}}_N).
\end{align}
Note that $\tilde{B}_0,\tilde{B}_r\in\R^{n\times n_\phi}$ and  $B_d  N_{uw}$ can be written as 
$
    B_d  N_{uw} = \tilde{B}_0 + \sum_{i=1}^n\sum_{j=1}^{n_\phi} \bm{e}_i^n \tilde{g}_{ij}  {(\bm{e}_j^{n_\phi})}^{\top},$ where $\tilde{g}_{ij}\in \mathbb{R}$ and $ |\tilde{g}_{ij}|\leq \tilde{B}_r (i,j).
$

The following result shows the feasibility equivalence of two LMIs that will be used for the robust stability of NNCS \eqref{close-sys}.

\begin{proposition} \label{prop:lmi_eq}
Consider the uncertain NNCS \eqref{close-sys} with $A_d, B_d$ satisfying the interval matrix uncertainty constraints \eqref{equ:uncertain}. Let the equilibrium point $\bm x_*=\bm 0$ and $(\bm x_*,\u_*,\bm v_*, \bm w_*)$ satisfy \eqref{equ:NNequ}, 
and let $\overline{\bm v}^1 \in \mathbb{R}^{n_1}$. Let $R_V$, $X(\bm \lambda)$, $\tilde{B}_0$ and $\tilde{B}_r$ be defined as in \eqref{rvrphi}, \eqref{eqxlambda}, and \eqref{equ:B_new_0}, respectively. Define $\hat n=2n+n_\phi$, $D = \left[A_r\; \tilde{B}_r\; \bm 0_{n\times n}\right]^\top$ 
and 
\begin{align}\label{equ:Zmat}
    Z(\bm \lambda, P) \!=\! \mat{\!-R_V^{\top} \!\mat{I_n \\ \bm 0_{m\times n}} \!P\! \mat{I_n & \bm 0_{n\times m}}\! R_V \!+\! X(\bm \lambda)& *\!\\ P\!\mat{A_0 & \tilde{B}_0} & -P\!}\!.
\end{align}
Then the following two statements are equivalent.
\begin{enumerate}[label=\alph*),wide]
    \item \label{it_LMI1} There exists a positive definite matrix $P\in \mathbb{S}^{n}_{\succ 0}$, a vector $\bm\lambda \in \mathbb{R}^{n_\phi}_{\geq 0}$, and positive scalars $\gamma_{ij}>0$ ($i\in \bbox{\hat n}, j \in \bbox{n}$) such that \mbox{\textup{(LMI-I)}} shown below is feasible:  
\begin{tcolorbox}[colback=white,top=0mm,bottom=-4mm,right=0mm,left=2mm]
\begin{flushleft}\mbox{\textup{(LMI-I)}} \end{flushleft} \vskip -2mm
\begin{align}
\label{equ:relax1-1}
& \mat{Z(\bm \lambda, P) \!+\!\sum\limits_{i=1}^{\hat n} \sum\limits_{j=1}^{n} \gamma_{ij} (D(i,j))^2 \bm e_i \bm e_i^\top & U \\ * &\! -V} \!\!\prec\! 0,\\ 
&\mbox{and \eqref{equ:nom2} holds}\nonumber
\end{align}\vskip-5mm
\end{tcolorbox}
where \vskip -5mm  \begin{subequations}\label{equ:VU}
    \begin{align}
        V &= \operatorname{diag}\left(\gamma_{11},\dots, \gamma_{1n},\dots, \gamma_{\hat n1},\dots,\gamma_{\hat nn}\right),\label{eqV}\\
        U &=\mat{\bm 0_{(n+n_\phi)\times n} & \cdots & \bm 0_{(n+n_\phi)\times n}\\ P & \cdots & P}\in\R^{\hat n\times \hat nn}.\label{eqU}
    \end{align}
    \end{subequations}\vskip -5mm
    \item \label{it_LMI2} There exists a positive definite matrix $P\in \mathbb{S}^{n}_{\succ 0}$, a vector $\bm\lambda \in \mathbb{R}^{n_\phi}_{\geq 0}$, and diagonal matrices $T\in \mathbb{R}^{\hat n\times \hat n}, S\in \mathbb{R}^{n\times n}$ such that \mbox{\textup{(LMI-II)}} shown below is feasible: 
\begin{tcolorbox}[colback=white,top=0mm,bottom=-4mm,right=0mm,left=2mm]
\begin{subequations}\label{equ:relax2}
\begin{flushleft}\mbox{\textup{(LMI-II)}} \end{flushleft}\vskip -2mm
\begin{align}\label{equ:relax2-1}
 & \mat{ \multicolumn{2}{c}{\multirow{2}{*}{$Z(\bm \lambda, P)+T$}} & * \\ \multicolumn{2}{c}{} & *\\ \bm 0_{n \times (n+n_\phi)} & P &  -S}\prec 0,  \\ \label{equ:relax2-2}
& D  S D^\top  \prec T,\\ \notag
& \mbox{and \eqref{equ:nom2} holds.}
\end{align}
\end{subequations} \vskip-5mm
\end{tcolorbox}
\end{enumerate}
\end{proposition}

\begin{IEEEproof}
(1) First, we prove b) $\Rightarrow$ a). 
Since \eqref{equ:nom2} is included in $\mbox{\textup{(LMI-II)}}$, we only need to show that given $P\in \mathbb{S}^{n}_{\succ 0}$, $\bm\lambda \in \mathbb{R}^{n_\phi}_{\geq 0}$, and diagonal matrices $T = \operatorname{diag}(t_1,\cdots,t_{\hat n})$, $S = \operatorname{diag}(s_1,\cdots,s_n)$ satisfying \eqref{equ:relax2-1}-\eqref{equ:relax2-2}, there exist positive scalars $\gamma_{ij}>0$ ($i\in \bbox{\hat n}, j \in \bbox{n}$) such that $P,\bm\lambda,\{\gamma_{ij}\}$ satisfy \eqref{equ:relax1-1}. 

This part of the proof will consist of three steps: i) we show that the diagonal elements of $T$ and $S$ are all positive; ii) we construct the candidate positive scalars $\{\gamma_{ij}\}$ using the diagonal matrices $T$ and $S$; iii) we prove \eqref{equ:relax1-1} holds with these constructed scalars $\{\gamma_{ij}\}$ and the same $P,\bm\lambda$.

$\mbox{i)}$ 
By the property of Schur complement of \eqref{equ:relax2-1}, we have $-S \prec 0$ and $$Z(\bm \lambda, P)+ T + \mat{\bm 0 & P}^\top S^{-1} \mat{\bm 0 & P} \prec0 .$$ Therefore, $S$ is positive definite, which implies that $s_1,\cdots,s_n>0$. From condition \eqref{equ:relax2-2}, we get $T\succ D S D^\top \succeq 0$, which implies that $T\succ 0$ and $t_1,\cdots,t_{\hat n} >0$.

$\mbox{ii)}$ Let $\epsilon'>0$ be an arbitrary positive real number, and denote matrix $$F\triangleq-\bm{1}_{n\times \hat n}SD^\top-DS\bm{1}_{n\times \hat n}^\top-\epsilon'\bm{1}_{n\times \hat n}S \bm{1}_{n\times \hat n}^\top.$$ Clearly, $F$ is symmetric, which implies all its eigenvalues are real. Since $T\succ D S D^\top$, $\lambda_{min}(T- D S D^\top)>0$. Thus, there exists a positive scalar $\epsilon$ satisfying $0<\epsilon<\epsilon'$, such that $$\lambda_{min}(T- D S D^\top)+\epsilon\lambda_{min}(F)>0.$$ 
By Weyl's inequality \cite[Theorem 4.3.1]{horn2012matrix}, we have $$\lambda_{min}(T- D S D^\top+\epsilon F)\geq \lambda_{min}(T- D S D^\top) + \epsilon \lambda_{min}(F) >0.$$ Therefore, the matrix $T- D S D^\top+\epsilon F$ is positive definite since it is symmetric and all of its eigenvalues are positive. Since $$\epsilon F-DSD^\top\prec -( D + \epsilon \bm{1}_{n\times \hat n}) S (D + \epsilon \bm{1}_{n\times \hat n})^\top,$$ we have $$T - ( D + \epsilon \bm{1}_{n\times \hat n}) S (D + \epsilon \bm{1}_{n\times \hat n})^\top\succ T-DSD^\top+\epsilon F\succ 0.$$ 
Denote $\tilde{D} \triangleq D + \epsilon \bm{1}_{n\times \hat n}$ and $G \triangleq T- \tilde{D} S \tilde{D}^\top$. Then, $G \succ 0$ and $\tilde{D} (i,j)> {D}(i,j) \geq 0$ for all $i\in \bbox{\hat n}, j \in \bbox{n}$.  By elementary matrix operations, we can get $$G(p,p) = t_p-\sum_{k=1}^n s_k (\tilde{D} (p,k))^2, \text{ for } p\in \bbox{\hat n},$$  and $$G(p,q) = -\sum_{k=1}^n s_k \tilde{D}(p,k) \tilde{D}(q,k), \text{ for } p,q\in \bbox{\hat n} \text{ and } p\not= q.$$
Define 
\begin{equation}
\gamma_{ij} = \frac{s_j \sum_{k=1}^{\hat n} \tilde{D} (k,j)C^G_{k,1} }{\tilde{D} (i,j) C^G_{i,1} },\;i\in \bbox{\hat n}, j \in \bbox{n}.\label{gamm}
\end{equation} 
Obviously, $\gamma_{ij}>0$ since $C^G_{i,1}>0$ from Lemma \ref{prop-pos} and $\tilde{D}(i,j)>0,s_j>0$, for $i\in \bbox{\hat n}, j \in \bbox{n}$. It is easy to check that for any $j \in \bbox{n}$, 
$$
\begin{aligned}
    \sum_{i=1}^{\hat n} \frac{1}{\gamma_{ij}} = & \sum_{i=1}^{\hat n} \frac{\tilde{D}(i,j) C^G_{i,1} }{s_j \sum_{k=1}^{\hat n} \tilde{D}(k,j)C^G_{k,1} } \\
    = &  \frac{\sum_{i=1}^{\hat n} \tilde{D}(i,j) C^G_{i,1} }{s_j \sum_{k=1}^{\hat n} \tilde{D}(k,j)C^G_{k,1} } =  \frac{1}{s_j}.
\end{aligned}
$$ 

When $i=1$, we have 
\begin{align*}
    & t_1 - \sum_{j=1}^n \gamma_{1j}(\tilde{D}(1,j))^2 >  t_1 - \sum_{j=1}^n \gamma_{1j} (\tilde{D}(1,j))^2 \\
    = & (\sum_{j=1}^{\hat n}  G(j,1) C^G_{j,1})/C_{1,1}^G  =  |G|/C_{1,1}^G > 0.
\end{align*}
When $i = 2,\dots,\hat n$, it holds that 
\begin{align*}
    & \sum_{j=1}^n \gamma_{ij} ({D}(i,j))^2 < \sum_{j=1}^n \gamma_{ij}(\tilde{D} (i,j))^2 \\
= & \frac{\sum_{j=1}^n s_j \tilde{D}(i,j)\sum_{k=1}^{\hat n} \tilde{D}(k,j)C^G_{k,1}}{C^G_{i,1}}\\
= & \frac{\sum_{j=1}^n \! s_j \tilde{D}(i,j)^2C^G_{i,1}\!\!+\!\!\sum_{k=1,k\not = i}^{\hat n}\sum_{j=1}^n \! s_j \tilde{D}(i,j)\tilde{D}(k,j)C^G_{k,1}}{C^G_{i,1}} \\
= & \frac{(t_i-G(i,i))C^G_{i,1}- \sum_{k=1,k\not = i}^{\hat n} G(k,i) C^G_{k,1}}{C^G_{i,1}} \\
= & \frac{t_i C^G_{i,1}-\sum_{k=1}^{\hat n} G(k,i)C^G_{k,1}}{C^G_{i,1}} = t_i,
\end{align*}
where the last equality is according to the Laplace expansion in Lemma \ref{lem:laplace}. 
Therefore, $\{\gamma_{ij}\}$ defined in \eqref{gamm} ensures that $\sum_{j=1}^n \gamma_{ij} (\tilde{D}(i,j))^2 \leq t_i$ and $\sum_{i=1}^{\hat n} \frac{1}{\gamma_{ij}} \leq \frac{1}{s_j}$, for $i\in \bbox{\hat n}, j \in \bbox{n}$. Thus,  
$$\sum_{i=1}^{\hat n} \sum_{j=1}^{n} \gamma_{ij}(D (i,j))^2 \bm e_i \bm e_i^\top \preceq T$$ and $$\diag(\sum_{i=1}^{\hat n}\frac{1}{\gamma_{i1}},\dots,\sum_{i=1}^{\hat n} \frac{1}{\gamma_{in}}) \preceq S^{-1}.$$

$\mbox{iii)}$ 
Note that 
\begin{align*}
    & Z(\bm \lambda, P) + \sum_{i=1}^{\hat n} \sum_{j=1}^{n} \gamma_{ij}(D (i,j))^2 \bm e_i \bm e_i^\top + UV^{-1}U^{\top} \\
    = & Z(\bm \lambda, P)  + \sum_{i=1}^{\hat n} \sum_{j=1}^{n} \gamma_{ij} (D (i,j))^2 \bm e_i \bm e_i^\top \\ 
      & + \mat{\bm 0 & P}^\top \diag(\sum_{i=1}^{\hat n}\frac{1}{\gamma_{i1}},\dots,\sum_{i=1}^{\hat n} \frac{1}{\gamma_{in}}) \mat{\bm 0 & P}\\
    \preceq & Z(\bm \lambda, P)\allowbreak + T + \mat{\bm 0 & P}^\top S^{-1} \mat{\bm 0 & P}  \prec  0.
\end{align*}
Using the Schur complement of $Z(\bm \lambda, P) + \sum_{i=1}^{\hat n} \sum_{j=1}^{n} \gamma_{ij}(D (i,j))^2 \bm e_i \bm e_i^\top + UV^{-1}U^{\top}$, it is clear that \eqref{equ:relax1-1} holds with such $P,\bm\lambda$ and $\{\gamma_{ij}\}$ given in \eqref{gamm}. This completes the first part of the proof.

(2) Next, we prove a) $\Rightarrow$ b). Similar to the first case, we only need to show that given $P\in \mathbb{S}^{n}_{\succ 0}$, $\bm\lambda \in \mathbb{R}^{n_\phi}_{\geq 0}$, and $\gamma_{ij}>0$ ($i\in \bbox{\hat n}, j \in \bbox{n}$) satisfying \eqref{equ:relax1-1}, there exist diagonal matrices $T$ and $S$ such that $P,\bm\lambda,T,S$ satisfy \eqref{equ:relax2-1}-\eqref{equ:relax2-2}. 

This part of the proof will consist of
two steps: i) we construct diagonal matrices $T^*$ and $S$ and show that they satisfy a relaxed version of \eqref{equ:relax2-1}-\eqref{equ:relax2-2}; ii) we prove that the feasibility of the relaxed conditions implies the feasibility of \eqref{equ:relax2-1}-\eqref{equ:relax2-2}.

$\mbox{i)}$ We first show that there exist diagonal matrices $T^*,S$ satisfying \eqref{equ:relax2-1} and $D S D^\top  \preceq T^*$. Let $V_i = \operatorname{diag} (\gamma_{i1},\gamma_{i2},\dots,\gamma_{in})$ for $i\in \bbox{\hat n}$. Then $V = \operatorname{diag}(V_1,V_2,\dots,V_{\hat n})$. Denote $D^{i}$ as a matrix whose entries are all zeros except that the $i$-th row of $D^{i}$ is the same as the $i$-th row of $ D$ for $i\in \bbox{\hat n}$. Obviously, $D = \sum_{i=1}^{\hat n} D^{i}$. 
Define diagonal matrices 
$$
\begin{aligned}
    T^* \triangleq & \operatorname{diag}(t_1,\dots,t_{\hat n}) = \sum_{i=1}^{\hat n} \sum_{j=1}^{n} \gamma_{ij} (D (i,j))^2 \bm e_i \bm e_i^\top, \\
    S \triangleq & \operatorname{diag}(s_1,s_2,\dots,s_n)
\end{aligned}
$$
where $s_j = \frac{1}{\sum_{i=1}^{\hat n} \frac{1}{\gamma_{ij}}}$, $j\in \bbox{n}$. It's easy to check that $T^* = \sum_{i=1}^{\hat n} D^{i} V_i (D^{i})^\top \succeq 0$ and $S = (\sum_{i=1}^{\hat n} V_i^{-1})^{-1} \succ 0$. 

Using the Schur complement of \eqref{equ:relax1-1}, we get
$$
\begin{aligned}
     & Z(\bm \lambda, P) + \sum_{i=1}^{\hat n} \sum_{j=1}^{n} \gamma_{ij} (D (i,j))^2 \bm e_i \bm e_i^\top + UV^{-1}U^{\top} \!\prec\! 0\\
\Rightarrow & Z(\bm \lambda, P)+ T^* + U \diag(V_1^{-1},\dots,V_n^{-1}) U^\top \allowbreak \prec 0.
\end{aligned}
$$
Since 
$$
\begin{aligned}
 U \diag(V_1^{-1},\dots,V_n^{-1}) U^\top 
 = & \mat{\bm 0 & P}^\top\! \sum_{i=1}^n V_i^{-1} \allowbreak \mat{\bm 0 & P}\\ 
 = & \mat{\bm 0 & P}^\top\! S^{-1}\! \mat{\bm 0 & P},
 \end{aligned}
$$ 
we have $$Z(\bm \lambda, P) + \mat{\bm 0 & P}^\top \!S^{-1}\! \mat{\bm 0 & P} \allowbreak + T^*\allowbreak \prec0. $$ Because $-S \prec 0$, the Schur complement of the aforementioned inequality implies that $S$ and $T^*$ satisfy \eqref{equ:relax2-1}, i.e., $$\mat{ \multicolumn{2}{c}{\multirow{2}{*}{$Z(\bm \lambda, P)+T^*$}} & * \\ \multicolumn{2}{c}{} & *\\ \bm 0 & P &  -S}\prec 0.$$ 

Next, we show $S$ and $T^*$ satisfy $D S D^\top \preceq T^*$. 
Let 
$$
H_i = \mat{D^{i} V_i (D^{i})^\top & D^{i}\\ * & V_i^{-1}} \text{ for } i\in \bbox{\hat n},\; H = \mat{T^* & D\\ * & S^{-1}}.
$$
Then, 
\begin{align*}
    H & = \mat{T^* & D\\ * & S^{-1}} = \mat{\sum_{i=1}^{\hat n}D^{i} V_i (D^{i})^\top & \sum_{i=1}^{\hat n}D^{i}\\ * & \sum_{i=1}^{\hat n} V_i^{-1}}  \\
     & =  \sum\limits_{i=1}^{\hat n} \mat{D^{i} V_i (D^{i})^\top  & D^{i}\\ * & V_i^{-1}}  = \sum\limits_{i=1}^{\hat n} H_i.
\end{align*}
For any $H_i$ where $i\in \bbox{\hat n}$, we have $V_i^{-1}\succ 0$ and the Schur complement of $H_i$ is $D^{i} V_i (D^{i})^\top - D^{i} (V_i^{-1})^{-1} (D^{i})^\top = \mathbf{0}$, which implies that $H_i \succeq 0$. Thus, we get $H  = \sum_{i=1}^{\hat n} H_i \succeq 0 $. By using the property of the Schur complement of $H$, we get $T^* - D S D^\top\succeq 0$.

$\mbox{ii)}$ 
Now we show that there exists a diagonal matrix $T$, together with $S$ above, satisfying \eqref{equ:relax2-1}-\eqref{equ:relax2-2}. Since $$Z(\bm \lambda, P)+ T^* + \mat{\bm 0 & P}^\top S^{-1} \mat{\bm 0 & P} \prec0,$$ it's easy to check that there always exists $\epsilon>0$ such that $$Z(\bm \lambda, P)+ T^* +\epsilon I_{\hat n} + \mat{\bm 0 & P}^\top  S^{-1} \mat{\bm 0 & P} \prec 0.$$ Let $T = T^* +\epsilon I_{\hat n}$. Then, we have $D S D^\top  \preceq T^* \prec T$ and 
$$ Z(\bm \lambda, P)+ T + \mat{\bm 0 & P}^\top S^{-1} \mat{\bm 0 & P} \prec 0 .$$ Using the Schur complement of the aforementioned inequality, $T$ and $S$ satisfy \eqref{equ:relax2-1} and \eqref{equ:relax2-2}. This completes the proof.
\end{IEEEproof}

The LMIs developed in \cite[Proposition 3.3]{ben2002tractable} can also be extended to certify the robust stability of the uncertain NNCS \eqref{close-sys}. Specifically, either statement in Proposition \ref{prop:lmi_eq} holds if and only if there exists $P\in \mathbb{S}^{n}_{\succ 0}$, $\bm\lambda \in \mathbb{R}^{n_\phi}_{\geq 0}$, $\gamma_{ij}>0$ ($i\in \bbox{\hat n}, j \in \bbox{n}$),  and $Y\in \mathbb{S}^{\hat n}$ such that \mbox{\textup{(LMI-III)}} shown below is feasible: 
\begin{tcolorbox}[colback=white,top=0mm,bottom=0mm,right=0mm,left=2mm]
\begin{subequations}\label{equ:relax3}
\begin{flushleft}\mbox{\textup{(LMI-III)}} \end{flushleft}\vskip -2mm
\begin{align}
\label{equ:relax3-1}
    & \mat{Y-\sum\limits_{i=1}^{\hat n} \sum\limits_{j=1}^{n} \gamma_{ij} (D (i,j))^2 \bm e_i \bm e_i^\top & U \!\\ * & V\!} \succ 0,\\ \label{equ:relax3-2}
    & Y\prec -Z(\bm \lambda, P), \\ 
    & \mbox{and \eqref{equ:nom2} holds},\nonumber
\end{align}
\end{subequations}
\end{tcolorbox}
\noindent 
where $Z(\bm \lambda, P)$ is defined in \eqref{equ:Zmat} and $U,V$ are defined in \eqref{equ:VU}. In fact, it's straightforward to verify the equivalence between $\mbox{\textup{(LMI-I)}}$ and $\mbox{\textup{(LMI-III)}}$ as \eqref{equ:relax1-1} is equivalent to \eqref{equ:relax3} from Schur complement, by substituting $Y$ with $-Z(\bm \lambda, P)$ in \eqref{equ:relax3-1}. 

Therefore, the three LMIs proposed above, $\mbox{\textup{(LMI-I)}}$, $\mbox{\textup{(LMI-II)}}$ and $\mbox{\textup{(LMI-III)}}$,  are equivalent. 
Furthermore, these three LMIs extend existing robust stability results for linear systems with state matrix uncertainties to NNCS with both state matrix and input matrix uncertainties: \eqref{equ:relax1-1} in $\mbox{\textup{(LMI-I)}}$ corresponds to the LMI in \cite[Theorem 1]{mao2003quadratic}, \eqref{equ:relax2} in $\mbox{\textup{(LMI-II)}}$ corresponds to the LMI in \cite[Theorem 4]{alamo2008new}, and \eqref{equ:relax3} in $\mbox{\textup{(LMI-III)}}$ corresponds to the LMI in \cite[Proposition 3.3]{ben2002tractable}. 
Additionally, the proof of feasibility equivalence shown above can be directly used to prove the feasibility equivalence of the three corresponding LMIs in \cite{alamo2008new,ben2002tractable,mao2003quadratic}. 
The relationship among these LMIs
is summarized in Fig. \ref{fig:rel}.

\begin{figure}[!t]
    \centering
    \includegraphics[width = 0.98\linewidth]{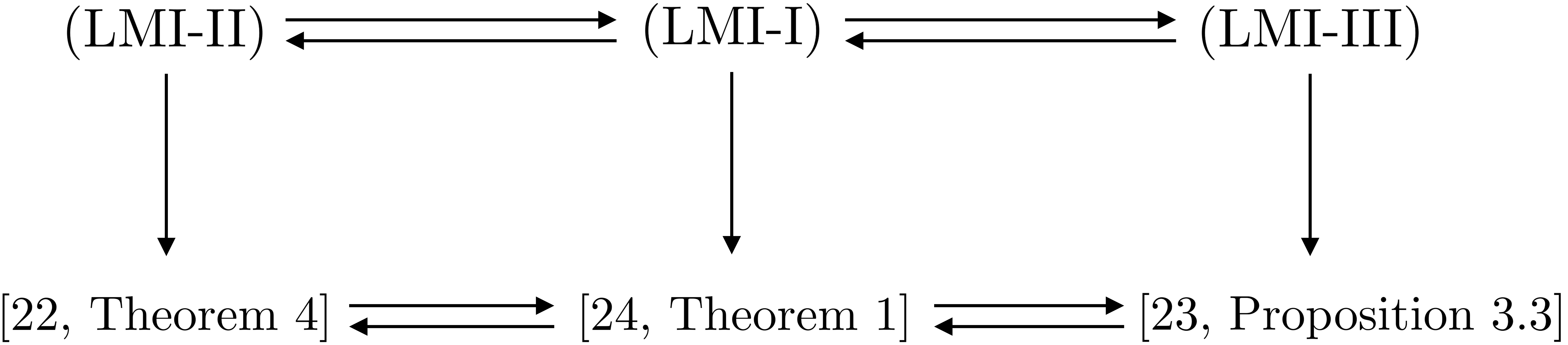}
    \caption{(LMI-I), (LMI-II), and (LMI-III) all serve as certificates for the robust stability of the NNCS in \eqref{close-sys} with interval uncertainties in both the  state and input matrices.  These LMIs extend the corresponding robust stability results for linear systems with interval state matrix uncertainties in  \cite{alamo2008new,ben2002tractable,mao2003quadratic}. The feasibility equivalence of the three proposed LMIs also implies the previously unestablished feasibility equivalence of the corresponding LMIs in \cite{alamo2008new,ben2002tractable,mao2003quadratic}.}
    \label{fig:rel}
\end{figure}

The following result shows that $\mbox{\textup{(LMI-I)}}$, $\mbox{\textup{(LMI-II)}}$ and $\mbox{\textup{(LMI-III)}}$ can be used to certify the robust stability of the uncertain NNCS \eqref{close-sys}.

\begin{theorem}\label{thm:relax1}
Consider the uncertain NNCS \eqref{close-sys} with $A_d, B_d$ satisfying the interval matrix uncertainty constraints \eqref{equ:uncertain}.
If either of the three LMIs - $\mbox{\textup{(LMI-I)}}$, $\mbox{\textup{(LMI-II)}}$ and $\mbox{\textup{(LMI-III)}}$ - is feasible, then the conclusion of Theorem \ref{thm:vertex} holds.
\end{theorem}

\begin{IEEEproof}
Since the feasibilities of the three LMIs are equivalent, we only need to show that one of the three conditions will lead to the robust stability of the uncertain system. For ease of readability, we choose $\mbox{\textup{(LMI-I)}}$. 

Suppose there exists a matrix $P \in \mathbb{S}_{\succ 0}^{n}$, a vector $\bm\lambda \in \mathbb{R}^{n_\phi}_{\geq 0}$, and positive scalars $\gamma_{ij}>0$ ($i\in \bbox{\hat n}, j \in \bbox{n}$) such that \eqref{equ:relax1-1} and \eqref{equ:nom2} of $\mbox{\textup{(LMI-I)}}$ hold. We will show that $P$ and $\bm \lambda$ also satisfy \eqref{equ:proplmi3} for any $A_d \in \vt([\underline{A}_d,\overline{A}_d])$ and $B_d \in \vt([\underline{B}_d,\overline{B}_d])$. Let $$\Xi = \mat{\sum_{i,j=1}^n \bm{e}_i^n  f_{ij} {(\bm{e}_j^n)}^{\top} \;&\; \sum_{i=1}^n\sum_{j=1}^{n_\phi} \bm{e}_i^n  \tilde{g}_{ij}  {(\bm{e}_j^{n_\phi})}^{\top}}.$$ 
    By Lemma \ref{lem:peterson} in Appendix, $\forall \gamma_{ij}>0 \;(i\in \bbox{\hat n}, j\in\bbox{n})$,
    \begin{align*}
        & \mat{-R_V^{\top} \mat{I_n & \bm 0}^{\top} P \mat{I_n & \bm 0}R_V + X(\bm \lambda)& * \\ P\mat{A_d & B_d}R_V & -P} \\
        = & Z(\bm \lambda, P) + \mat{\bm 0 & P}^\top \mat{\Xi & \bm 0} + \mat{\Xi & \bm 0}^\top \mat{\bm 0 & P}\\
        \preceq &  Z(\bm \lambda, P)+ \sum_{i=1}^{\hat n} \sum_{j=1}^{n} (\gamma_{ij}  (D (i,j))^2 \bm e_i \bm e_i^\top  \\
        & + \frac{1}{\gamma_{ij}} \mat{\bm 0 & P}^\top \bm e_j \bm e_j^\top \mat{\bm 0 & P})\prec 0.
    \end{align*}
    The last inequality is derived from the Schur complement of \eqref{equ:relax1-1}. So $P$ and $\bm \lambda$ satisfy \eqref{equ:proplmi3}, and thus \eqref{equ:nom1} by Proposition \ref{prop1}. Then the conclusion follows by Theorem \ref{thm:vertex}. 
\end{IEEEproof}

Theorem \ref{thm:relax1} indicates that the three LMIs - $\mbox{\textup{(LMI-I)}}$, $\mbox{\textup{(LMI-II)}}$ and $\mbox{\textup{(LMI-III)}}$ - are relaxations of $\mbox{\textup{(LMI-Vertex)}}$ with the same level of conservativeness.

\begin{remark}
    A brief comparison between \cite[Theorem 4]{alamo2008new} and \cite[Proposition 3.3]{ben2002tractable} was provided in \cite{alamo2008new}, which stated that \cite[Theorem 4]{alamo2008new} had the advantage of requiring less number of additional auxiliary decision variables. In this paper, we show that these two LMIs are actually equivalent in terms of feasibility. This equivalence indicates that \cite[Theorem 4]{alamo2008new} is more efficient than \cite[Proposition 3.3]{ben2002tractable} without sacrificing the feasibility. 
\end{remark}

\begin{remark}
The number of decision variables and the size of the matrices involved in solving the three LMI conditions are summarized in Table \ref{tab:complexity}. It can be observed that $\mbox{\textup{(LMI-II)}}$ has the least computation complexity in terms of the number of decision variables and the size of the LMIs when $n>1$. The computational complexity of $\mbox{\textup{(LMI-III)}}$ tends to be the largest, due to the extra decision variable $Y$. Despite the increased computational complexity, including additional decision variables, as in $\mbox{(LMI-I)}$ and $\mbox{(LMI-III)}$, can offer certain benefits depending on the application, such as improving the interpretability of intermediate variables and providing greater flexibility for further modification of the optimization problem.  For instance, the additional variable $Y$ in $\mbox{(LMI-III)}$ separates the constraints on the nominal system and the uncertainty bounds into two LMIs; by imposing additional constraints on $Y$, one can enforce desired properties  either on the nominal system or the uncertainties. The practical efficiency of the proposed LMI conditions is demonstrated through two numerical examples in Section \ref{sec:exp}.

\begin{table}[!t]
\caption{Comparison of $\mbox{\textup{(LMI-I)}}$, $\mbox{\textup{(LMI-II)}}$ and $\mbox{\textup{(LMI-III)}}$. The system state dimension is $n$, the number of neurons in the first layer is $n_1$, the total number of neurons is $n_\phi$ (i.e., $n_\phi = \sum_{i=1}^\ell n_i$), and $\hat{n} = 2n+n_\phi$.}
\label{tab:complexity}
\centering
\begin{tabular}{c|c|c}
\hline
 & \# Decision Variables & Size of LMIs \\ \hline

 $\mbox{\textup{(LMI-I)}}$ & $n(n\!+\!2\hat{n}\!+\!1)/{2}\!+\!n_\phi$ & $(\hat n\!+\!n_1)(n\!+\!1)$ \\ \hline

 $\mbox{\textup{(LMI-II)}}$   & $n(n\!+\!3)/{2}\!+\!\hat{n}\!+\!n_\phi$ & $n\!+\!2\hat{n}\!+\!n_1(n\!+\!1)$\\ \hline

 $\mbox{\textup{(LMI-III)}}$   & $(n\!+\!\hat{n})(n\!+\!\hat{n}\!+\!1)/{2}\!+\!n_\phi$ & $ (\hat n\!+\!n_1)(n\!+\!1) \!+\! \hat n$\\ \hline
\end{tabular}
\end{table}
\end{remark}

\begin{remark}
The problem of verifying the robust stability of uncertain NNCSs with an uncertain plant and an NN controller was also explored in \cite{yin2021stability}. However, the uncertainties discussed therein are categorized as structured uncertainties according to \cite{petersen2014robust}, characterized by IQCs, and the approach given in \cite{yin2021stability} is not directly applicable to the robust stability problem considered in this work. The uncertainties related to interval matrices in this work are less structured, only necessitating knowledge of their upper and lower bounds. 

On the other hand, in order to find the largest robust ROA inner-approximations, we can follow the same procedure as in \cite{yin2021stability} by adding $\operatorname{trace}(P)$ as the cost function of the LMIs developed before. For example, the optimization problem for $\mbox{\textup{(LMI-I)}}$ is formulated as follows:
\begin{equation}\label{equ:optim}
    \begin{aligned}
    \min_{P,\bm \lambda,\{\gamma_{ij}\}} \;\;\operatorname{trace}(P) \quad\text{s.t.}  \;\;\mbox{\textup{(LMI-I)}}\text{ holds}. 
\end{aligned}
\end{equation}
The optimization problem for other LMIs are formulated similarly. Based on Theorem \ref{thm:relax1} and the equivalence in feasibility of the three relaxed LMIs, the sizes of the robust ROA approximations are the same for $\mbox{\textup{(LMI-I)}}$, $\mbox{\textup{(LMI-II)}}$, and $\mbox{\textup{(LMI-III)}}$, which are smaller than the ROA approximation of $\mbox{\textup{(LMI-Vertex)}}$. However, the induced conservativeness of the relaxed LMIs is minimal compared to $\mbox{\textup{(LMI-Vertex)}}$, as demonstrated in the simulation examples.
\end{remark}

\begin{remark}
Although this work focuses on the robust stability analysis of uncertain NNCSs \eqref{close-sys} with linear dynamics and FNN controllers, the proposed method can also be extended to more general uncertain systems with nonlinear dynamics and other NN architectures (e.g., convolutional and recurrent neural networks), as QCs can describe a wide range of nonconvex and nonlinear functions used in general NNCSs \cite{megretski1997system,xu2020observer}. 
\end{remark}

\section{Simulation results} \label{sec:exp}
We use two simulation examples to illustrate the results of the preceding sections. In the following examples, the proposed LMIs are solved using MOSEK in MATLAB R2022b on a desktop with an Intel I7-8700K CPU and 32 GB memory.

\subsection{Inverted pendulum with uncertain length}
Consider the linearized inverted pendulum model: $$\dot{\bm{x}} = \mat{0 & 1\\ \frac{g}{\ell} & - \frac{\mu}{m \ell^2}} {\bm{x}} + \mat{0\\ \frac{1}{m \ell^2}} {\u}$$ 
where the state $\bm x =[\theta\;\dot{\theta}]^\top \in \mathcal{X} \subset \mathbb{R}^2$ represents the angular position $\theta$ and velocity $\dot{\theta}$, and $ \u \in \mathbb{R}$ is the control input. The state constraint set is given by $\mathcal{X}=[-2.5,2.5] \times[-6,6]$. The continuous-time dynamics is discretized with the sampling time $\Delta t=$ $0.02$ seconds. The inverted pendulum has the mass $m = 0.15$ kg and the friction coefficient $\mu = 0.05$ N$\cdot$m$\cdot$s/rad. We assume that there exists measurement uncertainty in the length of the pendulum such that  $\ell \in [0.5-\delta, 0.5+\delta]$, where $\delta$ is the level of uncertainty. The system thus contains interval matrix uncertainties that are in the form of \eqref{equ:uncertain}. We use the same NN controller as in \cite{yin2021stability} which was obtained through a reinforcement learning process using policy gradient. The NN controller is parameterized by a 2-layer FNN with 32 neurons per layer and \verb|tanh| as the activation function. The control input is saturated by $\u = sat^{\bar u}_{\underline u} \pi(\bm x)$ with $\bar u = - \underline u = 0.7$ N$\cdot$m. We also assume that $\bm v^1 \triangleq W^1\x+\bm b^1 \in [\underline{\bm v}^1,\overline{\bm v}^1]$ with $\overline{\bm v}^1 = -\underline{\bm v}^1 = 0.1\times \bm 1_{32\times 1}$.

\begin{figure}[!t]
\centering
  \includegraphics[width=0.85\linewidth]{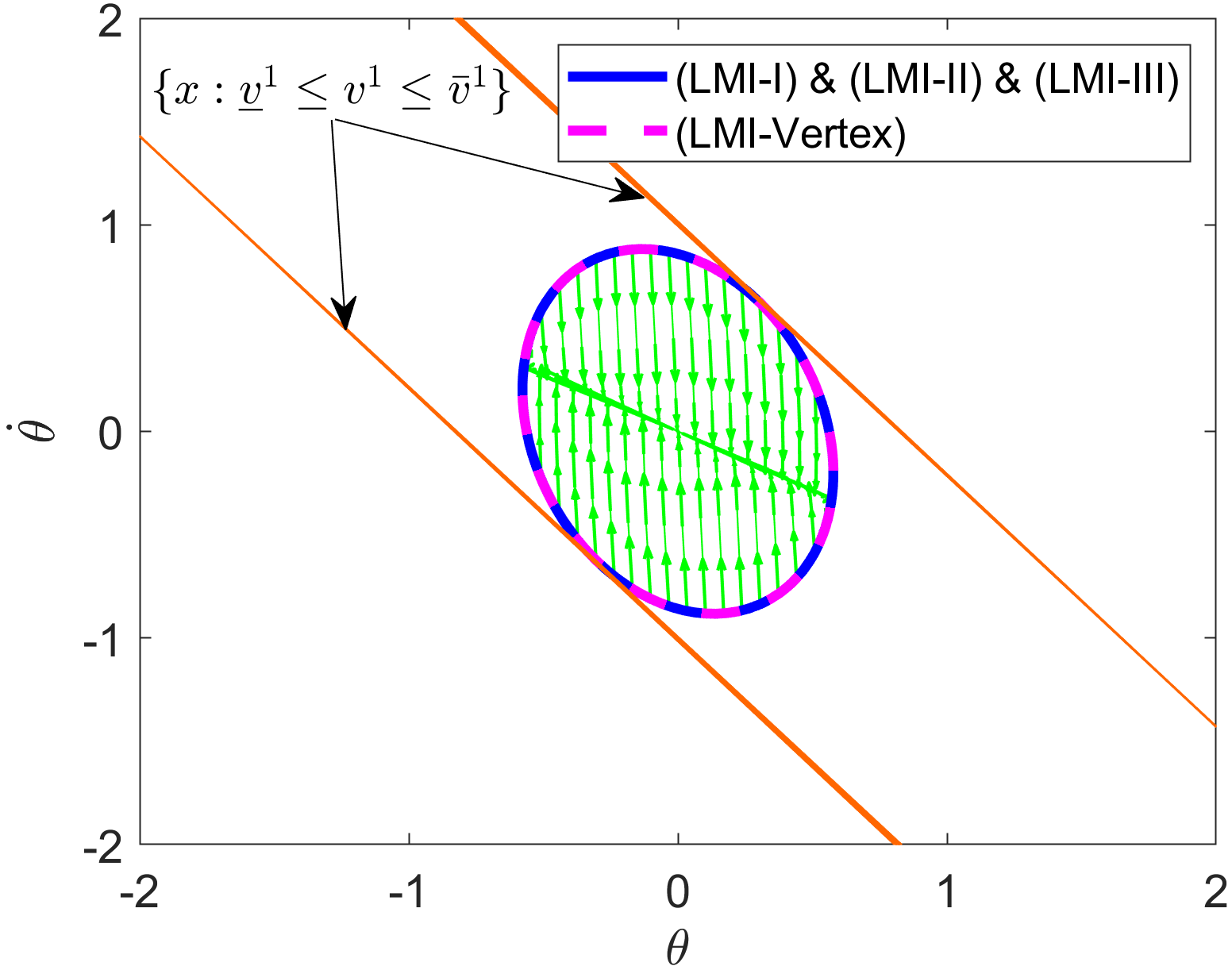}
  \caption{Inner-approximations of the robust ROA for the inverted pendulum model with uncertainty level $\delta = 0.01$. The three inner-approximations corresponding to the three LMIs are congruent. Trajectories with randomly selected initial states on the boundary are plotted in green. }
  \label{fig:exp1-1}
\end{figure}

\begin{figure}[!t]
\centering
  \includegraphics[width=0.85\linewidth]{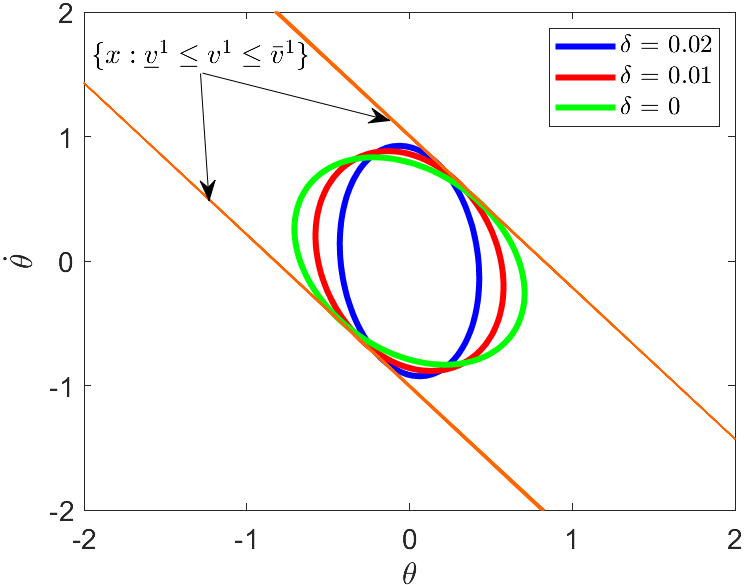}
  \caption{Inner-approximations of the robust ROA for the inverted pendulum model with varying uncertainty levels.
  }
  \label{fig:exp1-2}
\end{figure}

Fig. \ref{fig:exp1-1} shows inner-approximations of the robust ROA obtained by solving optimization problems \eqref{equ:optim} using the LMIs from Theorem \ref{thm:vertex} and Theorem \ref{thm:relax1}. 
The relaxed LMIs yield ellipsoid inner-approximations of equal size, within the set $\{\bm x: \underline{\bm v}^1 \leq \bm v^1 \leq \overline{\bm v}^1\}$ as enforced by \eqref{equ:nom2}, consistent with the LMI equivalence established in Section \ref{sec:relax}. The ROA approximations computed based on relaxed LMIs have a similar size to the one based on $\mbox{\textup{(LMI-Vertex)}}$ in Theorem \ref{thm:vertex}. Trajectories from randomly generated initial values on the inner-approximation boundary are plotted in green and all converge to the origin, consistent with Theorem \ref{thm:relax1}. 
The computation times for solving \mbox{\textup{(LMI-I)}}, \mbox{\textup{(LMI-II)}}, \mbox{\textup{(LMI-III)}}, and \mbox{\textup{(LMI-Vertex)}} are 1.886, 0.537, 4.301, and 2.991 seconds, respectively. The volumes of the ROA approximations are 1.560, 1.560, 1.560, and 1.561, respectively. 
Fig. \ref{fig:exp1-2} illustrates robust ROA inner-approximations under varying levels of uncertainty, revealing that the approximated ROAs shrink as the level of uncertainty increases.

\subsection{Mass-spring-damper with uncertain coefficients}

Consider the mass-spring-damper system consisting of $n_c$ carts as shown in Fig. \ref{fig:mass} \cite{jouret2023safety}.
Its state-space representation is 
$$
    \dot{\x}=\mat{\bm 0_{n_c \times n_c} & \bm I_{n_c} \\
-M^{-1} K & -M^{-1} C} \x+\mat
{\bm 0_{n_c \times n_c} \\
M^{-1}}\u
$$
where the state $\x = [z_1\;\cdots\;z_{n_c}\;\dot{z}_1\;\cdots\;\dot{z}_{n_c}]^\top\in\mathbb{R}^{2n_c}$ contains the position and velocity of the carts, the control input $\u = [f_1\;\cdots\;f_{n_c}]^\top\in\mathbb{R}^{n_c}$ combines the external forces applied on each cart, and $K,M,C$ are the same as those given in \cite{jouret2023safety}. The continuous-time dynamics is discretized with the sampling time $\Delta t=$ $0.1$ seconds. 
We assume that each cart has a mass of 1 kg, and the spring stiffness constants and damping coefficients are unknown and given in intervals: $k_i\in[1-\delta_k,1+\delta_k]$ and $c_i\in[0.1-\delta_c,0.1+\delta_c]$, $\forall i \in \li n_c \ri $ where $\delta_k=0.05$ and $\delta_c = 0.005$ reflect the level of model uncertainties. 

\begin{figure}[!b]
\centering
    \includegraphics[width=0.85\linewidth]{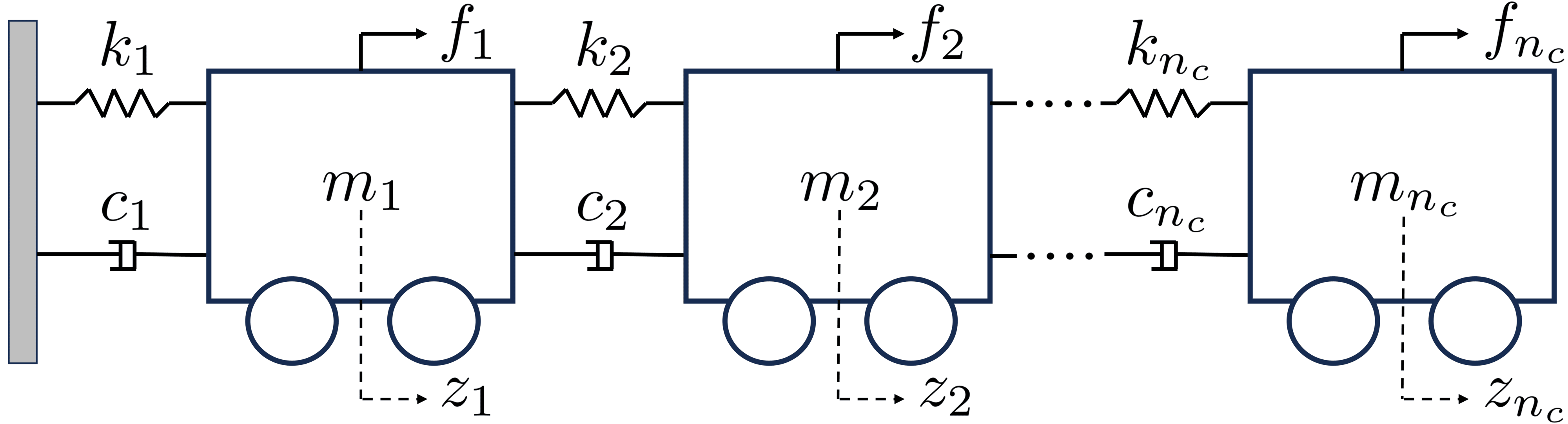}
\caption{Mass-spring-damper system with $n_c$ carts.}
\label{fig:mass}
\end{figure}

\begin{figure}[!t]
\centering
    \includegraphics[width=0.85\linewidth]{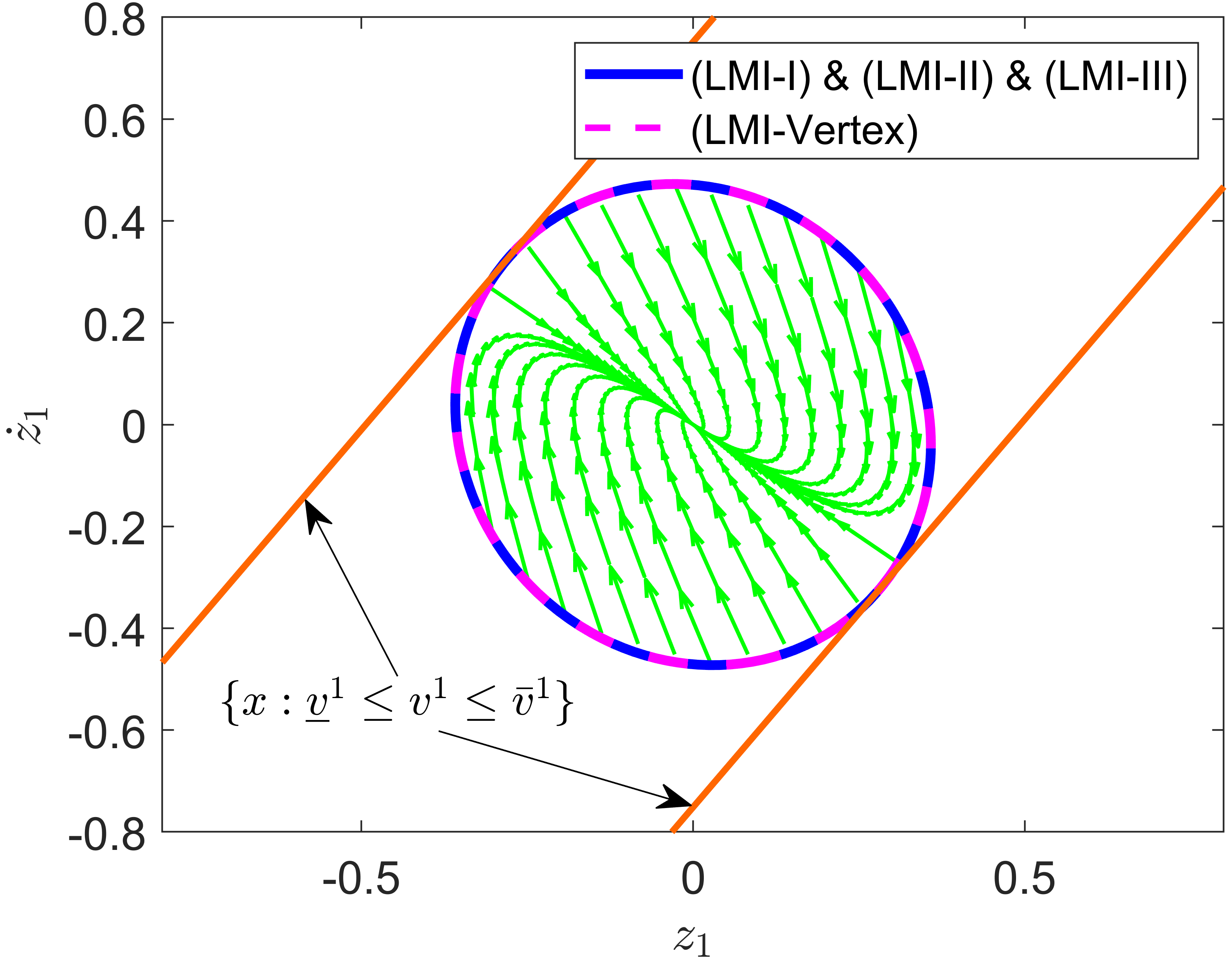}
\caption{Inner-approximations of the robust ROA for the mass-spring-damper model with $n_c = 1$. The three inner-approximations corresponding to the three LMIs coincide with each other. Trajectories with randomly selected initial states on the boundary are plotted in green.}
\label{fig:exp-mass}
\end{figure}

The FNN controller $\pi(\x)$ is trained using stochastic gradient descent to approximate an MPC controller that stabilizes the carts around the equilibrium point. It is parameterized by a 2-layer \verb|tanh|-activated FNN, with 8 neurons in each layer when $n_c = 1\text{ or }2$ and 16 neurons when $n_c = 3\text{ or }4$. It is assumed that $\bm v^1 \triangleq W^1\x+\bm b^1 \in [\underline{\bm v}^1,\overline{\bm v}^1]$ with $\overline{\bm v}^1 = -\underline{\bm v}^1 = 0.2\times \bm 1_{n_1\times 1}$. Fig. \ref{fig:exp-mass} shows the robust ROA inner-approximations on the ${z_1}-\dot{z}_1$ plane and the phase portrait of the closed-loop system in green, with $n_c = 1$. The results are consistent with the analysis in Section \ref{sec:relax}, as the three relaxed LMIs yield robust ROA inner-approximations with the same size, and all trajectories starting inside the robust ROA are driven to the origin despite the uncertainties. The comparison on runtimes and ROA approximation volumes for solving the three relaxed LMIs and the vertex-based LMI  with varying system dimensions are summarized in Table \ref{tab:massT} and \ref{tab:massV}. One can observe that the relaxed LMIs  achieve significant improvement in computational efficiency for high-dimensional systems with minimal conservativeness in ROA approximations compared to \mbox{(LMI-Vertex)}, which aligns with the results in Table \ref{tab:complexity}.

\section{Conclusion}\label{sec:concl}
In this paper, we investigated the robust stability problem for NNCSs with interval matrix uncertainties. Based on classic robust stability techniques and the QC-based descriptions of NNs, a novel LMI condition was proposed to certify the robust stability of uncertain NNCSs. Relaxed sufficient conditions based on LMIs were also presented which can reduce the computation burden involved in solving the LMIs. Feasibility of the three relaxed conditions was proved to be equivalent and their connections with existing robust stability results were also established.

\begin{table}[!t]
\centering
\caption{Comparison of runtime for solving the proposed LMI conditions for the mass-spring-damper example. The symbol ``-'' indicates that the computation did not finish within 3 hours.}
\label{tab:massT}
\centering
\begin{tabular}{c|cc|cc}
\hline
\multirow{2}{*}{} & \multicolumn{4}{c}{Time (s)}\\ \cline{2-5}

  & \multicolumn{1}{c|}{$n_c = 1$} & $n_c = 2$ & \multicolumn{1}{c|}{$n_c = 3$} & $n_c = 4$ \\ \hline

$\mbox{\textup{(LMI-I)}}$ & \multicolumn{1}{c|}{0.125} & 0.369 & \multicolumn{1}{c|}{2.692} & 6.332 \\ \hline

$\mbox{\textup{(LMI-II)}}$   & \multicolumn{1}{c|}{0.098} & 0.142 & \multicolumn{1}{c|}{0.169} & 0.320\\ \hline

$\mbox{\textup{(LMI-III)}}$   & \multicolumn{1}{c|}{0.152} & 0.409 & \multicolumn{1}{c|}{3.523} & 7.005 \\ \hline

$\mbox{\textup{(LMI-Vertex)}}$   & \multicolumn{1}{c|}{0.103} & 2.153 & \multicolumn{1}{c|}{673.4} & - \\ \hline
\end{tabular}
\end{table}

\begin{table}[!t]
\centering
\caption{Comparison of the sizes of the ROA approximations for the mass-spring-damper example.}
\label{tab:massV}
\centering
\begin{tabular}{c|cc|cc}
\hline
\multirow{2}{*}{} & \multicolumn{4}{c}{Volume of ROA approximation}\\ \cline{2-5}

  & \multicolumn{1}{c|}{$n_c = 1$} & $n_c = 2$ & \multicolumn{1}{c|}{$n_c = 3$} & $n_c = 4$ \\ \hline

$\mbox{\textup{(LMI-I)}}$ & \multicolumn{1}{c|}{0.531} & 288.2 & \multicolumn{1}{c|}{2036} & 22473 \\ \hline

$\mbox{\textup{(LMI-II)}}$   & \multicolumn{1}{c|}{0.531} & 288.2 & \multicolumn{1}{c|}{2036} & 22473\\ \hline

$\mbox{\textup{(LMI-III)}}$   & \multicolumn{1}{c|}{0.531} & 288.2 & \multicolumn{1}{c|}{2036} & 22473 \\ \hline

$\mbox{\textup{(LMI-Vertex)}}$   & \multicolumn{1}{c|}{0.531} & 288.9 & \multicolumn{1}{c|}{2051} & - \\ \hline
\end{tabular}
\end{table}

\bibliographystyle{IEEEtran}
\bibliography{final_ref}

\begin{thebibliography}{10}
\providecommand{\url}[1]{#1}
\csname url@samestyle\endcsname
\providecommand{\newblock}{\relax}
\providecommand{\bibinfo}[2]{#2}
\providecommand{\BIBentrySTDinterwordspacing}{\spaceskip=0pt\relax}
\providecommand{\BIBentryALTinterwordstretchfactor}{4}
\providecommand{\BIBentryALTinterwordspacing}{\spaceskip=\fontdimen2\font plus
\BIBentryALTinterwordstretchfactor\fontdimen3\font minus \fontdimen4\font\relax}
\providecommand{\BIBforeignlanguage}[2]{{%
\expandafter\ifx\csname l@#1\endcsname\relax
\typeout{** WARNING: IEEEtran.bst: No hyphenation pattern has been}%
\typeout{** loaded for the language `#1'. Using the pattern for}%
\typeout{** the default language instead.}%
\else
\language=\csname l@#1\endcsname
\fi
#2}}
\providecommand{\BIBdecl}{\relax}
\BIBdecl

\bibitem{fazlyab2020safety}
M.~Fazlyab, M.~Morari, and G.~J. Pappas, ``Safety verification and robustness analysis of neural networks via quadratic constraints and semidefinite programming,'' \emph{IEEE Transactions on Automatic Control}, vol.~67, no.~1, pp. 1--15, 2020.

\bibitem{hornik1989multilayer}
K.~Hornik, M.~Stinchcombe, and H.~White, ``Multilayer feedforward networks are universal approximators,'' \emph{Neural Networks}, vol.~2, no.~5, pp. 359--366, 1989.

\bibitem{nikolakopoulou2022dynamic}
A.~Nikolakopoulou, M.~S. Hong, and R.~D. Braatz, ``Dynamic state feedback controller and observer design for dynamic artificial neural network models,'' \emph{Automatica}, vol. 146, p. 110622, 2022.

\bibitem{schwan2023stability}
R.~Schwan, C.~N. Jones, and D.~Kuhn, ``Stability verification of neural network controllers using mixed-integer programming,'' \emph{IEEE Transactions on Automatic Control}, 2023.

\bibitem{zhang2022safety}
Y.~Zhang and X.~Xu, ``Safety verification of neural feedback systems based on constrained zonotopes,'' in \emph{IEEE Conference on Decision and Control}, 2022, pp. 2737--2744.

\bibitem{zhang2023backward}
Y.~Zhang, H.~Zhang, and X.~Xu, ``Backward reachability analysis of neural feedback systems using hybrid zonotopes,'' \emph{IEEE Control Systems Letters}, vol.~7, pp. 2779--2784, 2023.

\bibitem{zhang2024reachability}
------, ``Reachability analysis of neural network control systems with tunable accuracy and efficiency,'' \emph{IEEE Control Systems Letters}, vol.~8, pp. 1697--1702, 2024.

\bibitem{goodfellow2014explaining}
I.~J. Goodfellow, J.~Shlens, and C.~Szegedy, ``Explaining and harnessing adversarial examples,'' in \emph{International Conference on Learning Representations}, 2015.

\bibitem{fabiani2022reliably}
F.~Fabiani and P.~J. Goulart, ``Reliably-stabilizing piecewise-affine neural network controllers,'' \emph{IEEE Transactions on Automatic Control}, vol.~68, no.~9, pp. 5201 -- 5215, 2022.

\bibitem{karg2020stability}
B.~Karg and S.~Lucia, ``Stability and feasibility of neural network-based controllers via output range analysis,'' in \emph{IEEE Conference on Decision and Control}, 2020, pp. 4947--4954.

\bibitem{newton2023sparse}
M.~Newton and A.~Papachristodoulou, ``Sparse polynomial optimisation for neural network verification,'' \emph{Automatica}, vol. 157, p. 111233, 2023.

\bibitem{dawson2023safe}
C.~Dawson, S.~Gao, and C.~Fan, ``Safe control with learned certificates: A survey of neural {L}yapunov, barrier, and contraction methods for robotics and control,'' \emph{IEEE Transactions on Robotics}, vol.~39, no.~3, pp. 1749 -- 1767, 2023.

\bibitem{de2023event}
C.~de~Souza, A.~Girard, and S.~Tarbouriech, ``Event-triggered neural network control using quadratic constraints for perturbed systems,'' \emph{Automatica}, vol. 157, p. 111237, 2023.

\bibitem{hu2020reach}
H.~Hu, M.~Fazlyab, M.~Morari, and G.~J. Pappas, ``Reach-{SDP}: Reachability analysis of closed-loop systems with neural network controllers via semidefinite programming,'' in \emph{IEEE Conference on Decision and Control}, 2020, pp. 5929--5934.

\bibitem{jin2020stability}
M.~Jin and J.~Lavaei, ``Stability-certified reinforcement learning: A control-theoretic perspective,'' \emph{IEEE Access}, vol.~8, pp. 229\,086--229\,100, 2020.

\bibitem{pauli2021linear}
P.~Pauli, D.~Gramlich, J.~Berberich, and F.~Allg{\"o}wer, ``Linear systems with neural network nonlinearities: Improved stability analysis via acausal {Z}ames-{F}alb multipliers,'' in \emph{IEEE Conference on Decision and Control}, 2021, pp. 3611--3618.

\bibitem{yin2021stability}
H.~Yin, P.~Seiler, and M.~Arcak, ``Stability analysis using quadratic constraints for systems with neural network controllers,'' \emph{IEEE Transactions on Automatic Control}, vol.~67, no.~4, pp. 1980--1987, 2021.

\bibitem{yin2021imitation}
H.~Yin, P.~Seiler, M.~Jin, and M.~Arcak, ``Imitation learning with stability and safety guarantees,'' \emph{IEEE Control Systems Letters}, vol.~6, pp. 409--414, 2021.

\bibitem{petersen1987stabilization}
I.~R. Petersen, ``A stabilization algorithm for a class of uncertain linear systems,'' \emph{Systems \& Control Letters}, vol.~8, no.~4, pp. 351--357, 1987.

\bibitem{petersen1986riccati}
I.~R. Petersen and C.~V. Hollot, ``A {R}iccati equation approach to the stabilization of uncertain linear systems,'' \emph{Automatica}, vol.~22, no.~4, pp. 397--411, 1986.

\bibitem{petersen2014robust}
I.~R. Petersen and R.~Tempo, ``Robust control of uncertain systems: Classical results and recent developments,'' \emph{Automatica}, vol.~50, no.~5, pp. 1315--1335, 2014.

\bibitem{alamo2008new}
T.~Alamo, R.~Tempo, D.~R. Ram{\'\i}rez, and E.~F. Camacho, ``A new vertex result for robustness problems with interval matrix uncertainty,'' \emph{Systems \& Control Letters}, vol.~57, no.~6, pp. 474--481, 2008.

\bibitem{ben2002tractable}
A.~Ben-Tal and A.~Nemirovski, ``On tractable approximations of uncertain linear matrix inequalities affected by interval uncertainty,'' \emph{SIAM Journal on Optimization}, vol.~12, no.~3, pp. 811--833, 2002.

\bibitem{mao2003quadratic}
W.-J. Mao and J.~Chu, ``Quadratic stability and stabilization of dynamic interval systems,'' \emph{IEEE Transactions on Automatic Control}, vol.~48, no.~6, pp. 1007--1012, 2003.

\bibitem{kim2018standard}
K.-K.~K. Kim, E.~R. Patr{\'o}n, and R.~D. Braatz, ``Standard representation and unified stability analysis for dynamic artificial neural network models,'' \emph{Neural Networks}, vol.~98, pp. 251--262, 2018.

\bibitem{megretski1997system}
A.~Megretski and A.~Rantzer, ``System analysis via integral quadratic constraints,'' \emph{IEEE Transactions on Automatic Control}, vol.~42, no.~6, pp. 819--830, 1997.

\bibitem{xu2020observer}
X.~Xu, B.~A{\c{c}}{\i}kme{\c{s}}e, and M.~J. Corless, ``Observer-based controllers for incrementally quadratic nonlinear systems with disturbances,'' \emph{IEEE Transactions on Automatic Control}, vol.~66, no.~3, pp. 1129--1143, 2020.

\bibitem{gowal2019scalable}
S.~Gowal, K.~D. Dvijotham, R.~Stanforth, R.~Bunel, C.~Qin, J.~Uesato, R.~Arandjelovic, T.~Mann, and P.~Kohli, ``Scalable verified training for provably robust image classification,'' in \emph{Proceedings of the IEEE/CVF International Conference on Computer Vision}, 2019, pp. 4842--4851.

\bibitem{nemirovskii1993several}
A.~Nemirovskii, ``Several {NP}-hard problems arising in robust stability analysis,'' \emph{Mathematics of Control, Signals and Systems}, vol.~6, pp. 99--105, 1993.

\bibitem{calafiore2008reduced}
G.~Calafiore and F.~Dabbene, ``Reduced vertex set result for interval semidefinite optimization problems,'' \emph{Journal of Optimization Theory and Applications}, vol. 139, pp. 17--33, 2008.

\bibitem{jaulin2001interval}
L.~Jaulin, M.~Kieffer, O.~Didrit, and E.~Walter, \emph{Applied Interval Analysis}.\hskip 1em plus 0.5em minus 0.4em\relax Springer, 2006.

\bibitem{horn2012matrix}
R.~A. Horn and C.~R. Johnson, \emph{Matrix Analysis}.\hskip 1em plus 0.5em minus 0.4em\relax Cambridge University Press, 2012.

\bibitem{jouret2023safety}
L.~Jouret, A.~Saoud, and S.~Olaru, ``Safety verification of neural-network-based controllers: A set invariance approach,'' \emph{IEEE Control Systems Letters}, 2023.

\bibitem{meyer1978singular}
C.~Meyer~Jr and M.~Stadelmaier, ``Singular {M}-matrices and inverse positivity,'' \emph{Linear Algebra and Its Applications}, vol.~22, pp. 139--156, 1978.

\end{thebibliography}

\appendix

\begin{lemma}\label{prop-pos}
Given $A\in \mathbb{R}^{n\times n}$ as a positive definite matrix with strictly positive diagonal entries and strictly negative off-diagonal entries, then the cofactors of $A$ are strictly positive.
\end{lemma}

\begin{IEEEproof}
Since $A$ is positive definite, it is non-singular and its eigenvalues are positive. Thus, $A$ is a non-singular M-matrix \cite[Definition 1]{meyer1978singular}. Obviously, $A$ is irreducible since all entries of $A$ are non-zero and permutation does not introduce any zero entries. Since an irreducible non-singular M-matrix is strictly inverse-positive \cite[Theorem A.(ii)]{meyer1978singular},  $A^{-1}$ is entry-wise positive. 
Recall that $A^{-1} = \frac{adj(A)}{|A|}$, and the adjugate matrix of $A$ is the transpose of its cofactor matrix $C$  (i.e., the $(i,j)$-th entry of $adj(A)$ is $C_{i,j}^A$). 
Since $A$ is positive definite, $|A|>0$. Therefore, $adj(A) = |A|\cdot A^{-1}$ is entry-wise positive which indicates that all the cofactors are strictly positive. 
\end{IEEEproof}

\begin{lemma}\cite{petersen1987stabilization} \label{lem:peterson}
    Let $A$, $B$, $F$ be real matrices of suitable dimensions with $F^\top F \leq I$. Then, for any scalar $\gamma >0$, $AFB+B^\top F^\top A^\top \preceq \gamma A A^\top + \frac{1}{\gamma} B^\top B.$
\end{lemma}

\begin{lemma}[Laplace Expansion]\cite{horn2012matrix}\label{lem:laplace}
     Given a matrix $A \in \mathbb{R}^{n\times n}$, $|A| = \sum_{j=1}^n A(i,j) C^A_{i,j} = \sum_{j=1}^n A(j,i) C^A_{j,i}$ with $i\in \bbox{n}$, and $\sum_{j=1}^n A(k,j) C^A_{i,j} =  \sum_{j=1}^n A(j,k) C^A_{j,i} = 0$ with $i,k\in \bbox{n}, i\not=k$.
\end{lemma}

\end{document}